\documentclass{article}

\usepackage{amsthm, amssymb}
\usepackage{amsmath}
\usepackage{amsfonts}
\usepackage{epsfig,multicol}
\usepackage{pstricks}
\usepackage{graphicx}

\title{\bf A RESONANCE INTERACTION  OF \\
SEISMOGRAVITATIONAL MODES  ON TECTONIC PLATES.}

\author{Victor Flambaum$^{1,2}$, Gaven Martin$^{1}$ and Boris Pavlov$^{1}$ \thanks{Research supported in
part by grants from the Australian Research Council (VF) and the N.Z.
Marsden Fund (GM\&BP).  GM was supported in part by the Hausdorff Institute during the writing of this article.  \newline \newline
1. New Zealand Institute of  Advanced  Study,  Massey University,  Auckland, New Zealand. \newline
2. School of Physics, University of  New South Wales, Sydney, Australia   \newline\newline
The authors would like to thank the well informed referee who made number of helpful suggestions.
}  }

\date{\em In memory of our colleague and friend\\ Boris Sergeevich Pavlov}

\begin{document}

\maketitle
\newtheorem{theorem}{Theorem}[section]    
                                           
\newtheorem{lemma}[theorem]{Lemma}         
\newtheorem{corollary}[theorem]{Corollary} 
\newtheorem{remark}[theorem]{Remark}       
\newtheorem{definition}[theorem]{Definition}
\newtheorem{conjecture}[theorem]{Conjecture}
\newtheorem{proposition}[theorem]{Proposition}
\newtheorem{example}[theorem]{Example}

\newcommand{\param}{(\gamma,\beta,\beta')}
\newcommand{\parfour}{(\gamma,\beta,-4)}

\numberwithin{equation}{section}

\newcommand{\abs}[1]{\lvert#1\rvert}

\renewcommand{\theequation}{\thetheorem} 
                    
\makeatletter
\let \c@equation=\c@theorem

\begin{abstract}  This paper discusses resonance effects to advance a classical earthquake model,  namely the  celebrated  M8 global test algorithm.  This algorithmgives high confidence levels  for  prediction of  Time Intervals of Increased  Probability (TIP)
of an earthquake. It is based  on observation that  almost  80\%  of  earthquakes occur due to  the stress  accumulated from  previous earthquakes at the location  and   stored in form  of displacements against  gravity  and  static  elastic deformations
of the  plates. Nevertheless  the  M8 global test algorithm  fails to predict  some  powerful  earthquakes.  In   this paper  we suggest the additional  possibility of
considering  the dynamical  storage  of  the  elastic energy  on  the  tectonic plates  due  to  resonance   beats  of  seismo-gravitational
oscillations (SGO) modes of the plates.  We  make sure that  the tangential  compression  in  the middle plane of
an ``active zone''   of a tectonic plate may tune its  SGO modes  to  the resonance condition  of coincidence the frequencies
of the corresponding localized  modes with  the delocalized  SGO modes of the complement.
We  also consider the beats  arising between the modes under a  small  perturbations of the  plates, and, assuming  that the discord  between the
perturbed and unperturbed  resonance  modes  is  strongly dominated  by the  discord between  the non-resonance  modes estimate the energy transfer coefficient.  \end{abstract}

\maketitle

\newpage

\section{Exterior and  interior dynamics of tectonic plates.}

The  ``M8  Global test algorithm" (see \cite{KB} and references)  for earthquake prediction was designed  in 1984 at the International Institute of Earthquake Prediction
and Mathematical Geophysics (Moscow)  based  on the observation  that  almost 80 \% of  actual events  at the selected  location arise due to the stress
 built  up  thanks to  previous events at the   corresponding Earthquake-prone  (active)  zone. J. K. Gardner and L. Knopoff
observed a sequence of earthquakes in Southern California, with aftershocks removed,  was Poissonian, \cite{GK}.  Since then a mix of statistical and analytical techniques through the results of a global 20-year long experiment give indirect confirmation of  common features of both the predictability and the diverse behaviour of the Earth's naturally fractal lithosphere.
The statistics achieved to date prove (remarkably with confidence above 99 \%) the rather high efficiency of the M8 and M8-MSc predictions limited to intermediate-term middle- and narrow-range accuracy.  These models also adaptable,  see eg \cite{MMZK}.  There are other models,  we note those found in the works Turcotte and Schubert, \cite{TS} and Dahlen and Tromp \cite{DT} and our ideas may have implications for these as well,  though we do not give specific details here.
 
 \medskip
 
The analytical and mechanical  arguments used to derive the various models are based  on the assumption,
 that  for the most  part,  the  energy  at the active zone is  stored  in the form of   static  elastic deformation and  the displacement  of  tectonic plates in  the
 gravitational  field.
 
 \medskip

 Though both   M8 and the improved  MSc  algorithms are  extremely efficient  for prediction of the
   Time intervals of Increased Probability (TIP) of  earthquakes,  some  highly  dangerous events, such as the recent  Tohoku  earthquake (Japan,March  11,  2011)
    were  not  predicted.  In Tohoku
  the ``black box''  constructed based on the above  algorithms, removed the  TIP warning  from the list  of  expected earthquakes at the Tohoku location   70 days
   before the  earthquake, see the retrospective analysis  of the Global Test effectiveness given by  Kossobokov in \cite{Kossobokov_2013}.

The mechanical arguments for these  algorithms  are  derived  from  the idea  of quasi-statical (adiabatic) variation of the potential energy of the plates  during
 the   periods between earthquakes  while the tectonic plates  participate   in the ``exterior''  dynamics, such as floating of (fragments  of) the plates down the
 slopes formed  by  previous  earthquakes,   or  responding  to the  hydrodynamical  oscillations  in  the resonance  cavities   at the earthquake-prone (active) zones
   filled with magma, see for instance, \cite{SR_2004}.
   
   \medskip

On the other hand,  earthquakes  do arise  within a  background formed by  oscillations  of  the planet.  Many  seismologists  study these typical oscillations  with
  amplitudes  within $\sim  0.2 - 0.5 $ cm, and  periods  of circa  a few minutes.  However interesting  anomalies with periods  of circa  10  minutes   were noticed by  G.A.Sobolev
   and A.A.Lyubushkin, see \cite{Sobolev}, in the course of their analysis  of seismological data preceding  the Sumatra earthquake on  December   26, 2004, recorded {\it in the
    remote zone}  of the earthquake.  Moreover,  some  decaying periodic  patterns, see Fig. 4  in \cite{Sobolev}, were noticed on the relevant spectral-temporal
     (time-spectral) cards, see below, with  the periods $\sim 100$  min. On  diagram  7  of the same paper one can see decaying  patterns with even greater
     periods  $\sim 2 400$  min  arising prior  the  major earthquake. Unfortunately the  authors of  \cite{Sobolev} did not consider  these as  important  details and  did  not
     develop  any extended  analysis  of   them in their paper.
     
     \medskip

Approximately a decade before this, see \cite{LPO90},  long periodic oscillation  patterns,  with periods  $\sim  40-70$  minutes, were registered  by  E.M. Linkov
 using a ``vertical pendulum'', constructed  especially  for  studying  Seismogravitational  Oscillations (SGO) of the Earth.  These oscillations  have been intensely
  monitored  during  the last decade,  see for instance \cite{P2002,Petrova_FZ08}, as  an  important  component of  the ``interior dynamic'' of the plates. They form
   a  natural dynamical background of  catastrophic events such as   Earthquakes, Tsunami and  Volcanic Eruptions.  Monitoring of  SGO confirmed the hypothesis \cite{Thesis,PP2008}
of their  spectral nature.  According to   \cite{PP2008},   the  SGO  should be interpreted   as decaying  flexural (vertical)  eigen-modes  of large tectonic plates  with
    linear size  up to few  thousands kilometers.

    The typical energy  of  the  mode  may be estimated  based  on spectral (frequency),  physical properties such as density and Young's modulus
    and  geometric characteristics of the plate.  For instance, the elastic energy stored in a single  SGO mode with frequency 200 $\mu$HHz and amplitude   $2\times 10^{-3}$ m
     on a tectonic plate with area circa  $ 10^{14}$ m$^2$, thickness 10$^5$ m and density 3380 kg m$^{-3}$ is estimated as $54\times 10^{9}$ joules.  This is almost equivalent
      to the seismic moment (``full  energy'') of the 4M  earthquake in Johannesburg (South Africa) November  18, 2013.
      
      \medskip

        While discussing the inner dynamics  of tectonic plates
       in our recent  publications,   see  \cite{MP_2014,FP_arxiv_2015} and similarly \cite{IMPP_2012},  we  proposed to  take  into account  the migration of elastic energy  between regions of tectonic plates, caused  by
        beating of the resonance  spectral modes  localized  on the regions.  For an active  zone, already  unstable  under  statical stress, the  migration  of energy
        defined by the  resonance  beating  might  be sufficient to  trigger an earthquake.

        We therefore suggest that the modelling  of  tectonic processes  with regard  of  resonance
         migration of energy would probably help in developing more  realistic  theoretical scenarios for an earthquake. In the simplest case  of  two resonance  SGO  modes, localized
          on  neighbouring  regions $\Omega_{\varepsilon},\, \Omega_c$, the  beating pattern is periodic and the  amount  of elastic energy transferred from one location to
          another on each period  is defined by the corresponding {\em transfer coefficient.}  The  transfer coefficient  may be large
for exact tuning of the corresponding  frequencies  $|\nu_{\varepsilon} -\nu_c | << \bar{\nu}\equiv |\nu_{\varepsilon} + \nu_c |/2$. Thus we study the influence of
resonance effects for this geophysical phenomenon and how t can be used to inform advances in the classical model and give examples which demonstrate that resonances can be a reason for
earthquake. However we are not able to present a full model for the phenomenon.  Some discussions pertain to
one-dimensional systems and cannot present a completely realistic model, but they demonstrate
the possibility of such mechanism in seismic phenomena as quite natural values of parameters are used. We hope the interested reader can follow the construction of the model and come to understanding
of how full implementation might be realised. The examples we offer can be considered as benchmarks for
future,  more realistic models,  based on these suggested ideas.

\medskip 
 In  this  paper  we  aim at
estimating   the  transfer coefficient, while considering  the  problem in frames of  perturbation analysis, depending  on the spectral characteristics  of  the  unperturbed
 modes  on disjoint  regions and  the type of  interaction imposed.   Usually  the  SGO  modes  are  presented  on  time-spectral  cards  obtained from the   corresponding  seismograms  via averaging  of  the   SGO  amplitudes, with  certain
  frequency,  on a step-wise  system of  time - windows, obtained  by  shifting  an initial window  by certain  interval of time  on each  step. The boundaries of the
  spectral -time  domains  on the cards, where the averaged  amplitude  of the  SGO mode  exceeds the  given  value $A$, form a system of isolines in the frequency/time
   coordinates $\nu,T$.  The horizontal axis for time, is  graded in hours, the vertical axis, for the frequencies, is  graded in $\mu$Hz.   Doctor  L. Petrova  provided us
   with  some  time-spectral cards from her private collection  and  shared  with us  some useful and interesting  comments
concerning the interpretation  of  the  cards  in terms  of SGO  dynamics, which  we referred  in our previous  publication  \cite{FP_arxiv_2015}. Find below  one of
 her cards which was obtained   by her from  seismograms  recorded on  SSB  station, France, during  the period  preceding  the  powerful earthquake
  on 26  September  2004  in Peru.
 The   averaging  of  the  amplitudes of the  seismo-gravitational oscillations, with  certain frequency,
was  done, after appropriate  filtration, on a  system of 20 hours  time - windows, obtained  by  shifting  an initial window  by  30 minutes on each  step. The  relief
of the  window-averaged  squared amplitude on the cards  is graded  by the isolines, with the step  $\delta A^2 =  \frac{1}{10} \left[ A^2_{max}  - A^2_{\min}\right]$, and
 is painted  accordingly  between the isolines with  shades depending  on  the square amplitude  $A^2$ : dull grey for the  background  value  $A^2_{\min}$   and white for
  the maximal value  $ A^2_{max}$. See  more comments in \cite{FP_arxiv_2015}.
\begin{center}
 \includegraphics [width = 3 in] {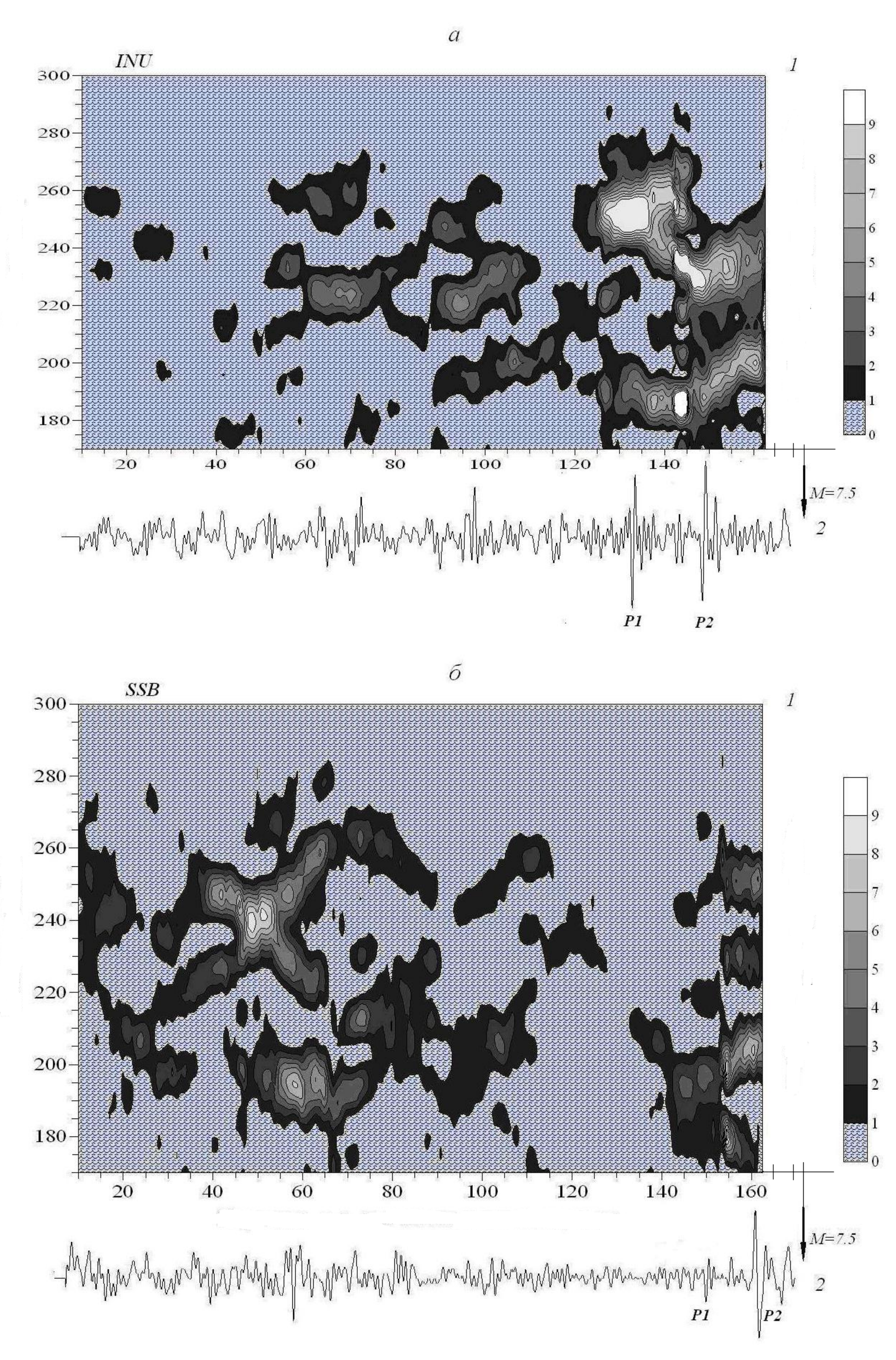}\\
  { \bf A  time - spectral card }
\end{center} 

The most interesting kind of SGO was represented  by  {\it pulsations}, which were observed by
    Linkov's team as intense  short (4-6  hours)  and  sometimes  repeated  in 30-100 hours, pulses of SGO with large  amplitudes.  Pulsations have been  registered  before 95 \%  of powerful earthquakes and  may be  considered  as natural precursors of them, see \cite{LPO90}.   One may hope that a  deeper analysis
      of the microseismic data  in \cite{Sobolev} would  reveal a connection  between the above  ``anomalies'' with SGO  and  pulsation  patterns studied  by Linkov et al.
    Petrova  attracted  pointed out to us  a peculiar  detail on the  above SSB card above, consisting of   two  groups  of stationary modes with  almost
       equal frequencies and visually similar relief   in $ \Delta^{SSB}_2 = (190,200)\times(55,65)$ and  in $\Delta^{SSB}_3 = (200,210)\times(145,165)$. The  pair was
         interpreted by    by   her  as a typical  ``seismo-gravitational  pulsation''.
         
         \medskip

In \cite{FP_arxiv_2015}  pulsations are interpreted as beatings of spectral modes  on the  tectonic plates, arising
due to  resonance interaction  of the  SGO modes, with close frequencies, while  some of them  are localized   on active zone  $\Omega_{\varepsilon}$  and   others - on
 the complement $\Omega_c$. According to classical  mechanics, \cite{Landau_Mechanics}, the  resonance between two ``oscillators'' $\Omega_{\varepsilon, \Omega_c}$  with
  close frequencies  $\omega_{\varepsilon},\,\omega_c$ and  precise tuning, the beating  of a pair  of modes  defines  the periodic  energy  migration between the  regions
   $\Omega_{\varepsilon, \Omega_c}$: $E_{\varepsilon} (t) =  \bar{E}_{\varepsilon} + \delta E (t),\, E_{c} (t) =  \bar{E}_{c} + \delta E (t)$.  Here the migrating part
    $\delta E (t)$ of the total energy   arises, with opposite phases, in  both locations, with  total energy conserved $E = \bar{E}_{\varepsilon} + \bar{E}_{c}$.

 While two  weakly connected  oscillators in resonance   yield a periodic beating,  a  larger  group  of    oscillators, under  resonance conditions
 $ |\omega _i - \omega_k |\equiv \delta_{ik} << \overline{\Omega} $
with respect to the  average  frequency $\bar{\omega} \equiv \frac{1}{N}\sum_{k=1}^N \,\, \omega_k $
reveal  a  quasi-chaotic  beating  phenomenon, if  the  difference frequencies $\delta_{ik}$  are non-co-measurable.

The problem  of  estimation of  the energy transfer  associated  with beats in the system of several  connected oscillators can be  reduced to  the similar problem
for a single oscillator  under almost  periodic resonance force.  This  is  treated  in Landau's book \cite{Landau_Mechanics}, Section 22. Indeed, we will consider  two  1D  oscillators,
with  masses  $m, M$ attached  to  springs $v, V$ and  connected  by an hermitian pair of  elastic bonds  $\gamma, \Gamma$ constrained by  the Hermitian  requirement
$\gamma \,\,m^{-1} =  \Gamma \,\, M^{-1} = \varepsilon $. For instance,  a pair of oscillators the dynamics is described by  the equations:
\[
x^{''} +  \frac{v}{m} x  + \gamma\,\, m^{-1} X = 0,
\]
\begin{equation}
\label{pairwise_beating}
X^{''} + \frac{V}{M} X + \Gamma\,\, M^{-1} x = 0.
\end{equation}
While the elastic bonds  are neglected, the eigenfrequencies  of the oscillators can be calculated  to be $\omega_m = \sqrt{v \,m^{-1}} \equiv \sqrt{\lambda_m},\,\,
\omega_M =  \sqrt{V \,M^{-1}} \equiv \sqrt{\lambda_M }$,  but with regard of the bonds are found  from  the quadratic equation  for  $\lambda = \omega^2$
\begin{equation}
\lambda_{\pm} =  \frac{\lambda_m + \lambda_M}{2}\pm \sqrt{ [\frac{\lambda_m + \lambda_M}{2}]^2  -   \lambda_m \lambda_M  + \varepsilon^2 } =
\frac{\lambda_m + \lambda_M}{2}\pm \sqrt{ [\frac{\lambda_m - \lambda_M}{2}]^2  + \varepsilon^2 }
\end{equation}
Hereafter we  will assume that {\em the bonds are  relatively weak} so that  we may  calculate the frequencies $\omega_{\pm} \equiv \sqrt{\lambda_{\pm}}$ of the normal modes
of the pair approximately  based on $\varepsilon << [\frac{\lambda_m - \lambda_M}{2}] \equiv \delta > 0 $, at least up to second order  with respect  to
$\varepsilon^2 [2 \delta]^{-1}$:
\begin{equation}
\label{approx_lambda}
\lambda_{\pm} = \frac{\lambda_m + \lambda_M}{2} \pm \frac{\lambda_m - \lambda_M}{2} \pm  \frac{\varepsilon^2}{2 \delta}.
\end{equation}
This implies, up to second order  with respect  to  $\varepsilon^2 [2 \delta]^{-1}$:
\[
\lambda_+ \approx \omega^2_m + \varepsilon^2 \, [2\delta]^{-1},\,\, \lambda_- = \omega^2_M -  \varepsilon^2 [2\delta]^{-1}
\]
and allows to calculate the complex  normal modes with  $\omega_{\pm} = \sqrt{\lambda_{\pm}}$
\begin{equation}
\label{normal_modes}
\left(
\begin{array}{cc}
e_+ \\
E_+
\end{array}
\right) \, e^{\pm i \omega_+ t},\,\, \left(
\begin{array}{cc}
e_- \\
E_-
\end{array}
\right) \, e^{\pm i \omega_- t}.
\end{equation}
The  eigenvectors $\left(
\begin{array}{cc}
e_+ \\
E_+
\end{array}
\right) $ and $\left(
\begin{array}{cc}
e_- \\
E_-
\end{array}
\right) $
are found from the homogeneous equations
\begin{equation}
\label{Eq : eigenvectors_pm}
\left(
\begin{array}{cc}
- \frac{\varepsilon^2}{2\delta} & \varepsilon\\
\varepsilon & - \delta - \frac{\varepsilon^2}{2\delta}
\end{array}
\right)\,\, \left(
\begin{array}{c}
e_+\\
E_+
\end{array}
\right) =  0,\,\,\,
\left(
\begin{array}{cc}
\delta + \frac{\varepsilon^2}{2\delta} & \varepsilon\\
\varepsilon &  \frac{\varepsilon^2}{2\delta}
\end{array}
\right)\,\, \left(
\begin{array}{c}
e_-\\
E_-
\end{array}
\right) =  0,
\end{equation}
which yield for the normalized  eigenvectors, up  to  $ O (\varepsilon^2 \,\,/\delta^2)$:
\begin{equation}
\left(
\begin{array}{c}
e_+\\
E_+
\end{array}
\right) = \left(
\begin{array}{c}
1\\
\frac{\varepsilon}{2\delta}
\end{array}
\right),\,\,
\left(
\begin{array}{c}
e_-\\
E_-
\end{array}
\right) =
\left(
\begin{array}{c}
-\frac{\varepsilon}{2\delta}\\
1
\end{array}
\right).
\end{equation}

In \cite{Landau_Mechanics}  the beats phenomenon  is  considered  for  a harmonic  external force. Hereafter we consider above  system of  two oscillators $x, X$, governed
by the  equations  (\ref{pairwise_beating}), under  initial condition  $x(0) = 1, x'(0) = 0, X(0) = \frac{\varepsilon}{2\delta}, X'(0) = 0$. These  initial conditions  correspond to  an excitation
of  the   normal mode  of above  the pair oscillators, with second order terms $\varepsilon^2 \delta^{-2}$ neglected:
\begin{equation}
\label{normal_mode_+}
\left(
\begin{array}{c}
x_+ (t)\\
X_+ (t)
\end{array}
\right) = \left(
\begin{array}{c}
1\\
\frac{\varepsilon}{2\delta}
\end{array}
\right) \cos \omega_+ t.
\end{equation}
Then  the dynamics if the first oscillators  may be considered  under  an exterior  force :
\begin{equation}
\label{exterior_force}
x^{''} + \omega^2_m x  +  \frac{\varepsilon^2}{2\delta} \cos \omega_+ \, t = 0,\,\, \mbox{with}\,\,\, \omega_+ = \omega_m  + \frac{\varepsilon^2}{2 \delta },
\end{equation}
under  the  initial conditions $x (0) = 1, x' (0) = 0.$  Following  \cite{Landau_Mechanics} for the real solution of the above equation  (\ref{normal_mode_+}) we introduce 
a  complex characteristic
\begin{equation}
\label{characteristic_1}
\xi = x'  +  i \omega_m x,
\end{equation}
and  rewrite the above equation (\ref{exterior_force}) for $\xi$ as
\begin{equation}
\label{characteristic_2}
\frac{d\xi}{dt} - i\omega_m  +  \frac{\varepsilon^2}{2\delta} \cos \omega_+ \, t,
\end{equation}
and  obtain the  solution $\xi (t)$ as  in  \cite{Landau_Mechanics}:
\begin{equation}
\xi (t) =  e^{i\omega_m t}\left[ \int_0^t  \frac{\varepsilon^2}{2\delta} \cos \omega_+ \, \tau \,\,\, e^{i\omega_m \tau}  d\tau  +  i  \omega_m \right]
\end{equation}
The energy ${\mathcal{E}} (x)$ of the oscillator $x$ is can be calculated as per \cite{Landau_Mechanics} again,   in terms  of this  characteristic as
\begin{equation}
\label{energy_x}
\frac{m}{2}  |\xi (t)|^2 =  \frac{m}{2}  \left[ |x'|^2 + \omega_m^2  |x|^2 \right] = {\mathcal{E}} (x).
\end{equation}

For two  connected oscillators $\delta E (t)$,  the migrating  part  of the total energy, is  close  to  zero  if the  difference in frequency
is relatively  large. Vice versa, it may be  close  to the full energy  if the tuning  is  very  sharp, that is  difference frequency is relatively small,
$|\omega_i - \omega_k| \equiv \delta_{ik} <<  \overline{\omega}$.  An example of this is  the  celebrated Wilberforce pendulum, \cite{L.Wilberforce, Wilberforce,Pendulum}.  We  posit that, in  the case of  tectonic plates,  beating  of the
  resonance SGO modes   implies (depending on conditions) migration of  an essential  part  $\max \delta E (t) = k \,\, E$ of the total energy $E$, with the
   transfer coefficient $ k = k(\delta \omega/\omega)$, depending on the relative  variation  of the frequencies. Then  the  total  energy   of the active zone
    $\Omega_{\varepsilon}$  at some moment  $T_0$ in time may exceed the destruction  limit  $E_{\varepsilon} (T) =  \bar{E}_{\varepsilon} + \delta E (T_0) > E_{\varepsilon}^d$,
    causing  the destruction of  some structure in   the active zone and  thus  triggering an earthquake.
    
    \medskip

Exactly this  scenario was considered in \cite{FP_arxiv_2015}  as  a resonance
     mechanism for an earthquake. Comparison of the full energy  of a single  SGO mode with the seismic moment requires an  estimation  of the  transfer  coefficient   $k$. This   becomes an  important  question  in analytic modelling  of  the resonance  mechanism  for an earthquake.
     
To obtain this  estimation of the transfer coefficient  $k$, in  the next Section  2  we  sketch the basics  of the
 Kirchhoff model  of the thin  plate, which is used  hereafter as the  tectonic plate. We leave the matter of the hydrodynamical component  of the  dynamics and the dissipation  to future work,
and concentrate  our  attention here on  the  effect  of the  compressing (tangential) tension  in the middle plane of the  plate $\Omega_{\varepsilon}$, causing a
  lowering  of the  eigenfrequencies of the active zone  to the  resonances with  the  SGO  modes of the complement.  Then in  Section 3 we  consider  the resonance condition  for
   the circular  active zone,  with the  prescribed  frequency  $200 \mu Hz$, choosing typical physical and geometrical parameters  of a pair $\Omega_{\varepsilon},\,\Omega_c$
   of the disjoint  plates. We consider the  simplest  explicitly solvable  example  of   two  unperturbed circular disjoint  plates,
and make sure that the resonance condition is fulfilled for the  active zone  $\Omega_{\varepsilon}$,
   under the compressing  tension,  with the parameters properly chosen, and a  circular plate  $\Omega'_c$, while the  tangential compression on the large plate is neglected
   in  the  resonance  interval  of frequencies. The circular large plate $\Omega'_c$ in Section 3  will differ from  the actual  complement $\Omega_c$  of the  circular active zone $\Omega_{\varepsilon}$  in minor ways,
 dominated  by  the typical wavelengths of  lower  SGO modes  on the plate. Hence so are the resonance conditions $\omega'_c\approx 200 \mu Hz \approx \omega_c$. A more accurate
  analysis  of the ring-like
complement  will be postponed to   an appendix, where  the Neumann-to-Dirichlet map  of the  ring is  calculated  in terms of Bessel functions. Similar  arguments will work for the
sectorial  boundary active  zone.

\vskip0.3cm
\section{Modelling the resonance  interaction
between SGO   and   beating  phenomena. }

This article investigates theoretical aspects  of  the  resonance  interaction of  SGO modes  based  on  the  Kirchhoff  model  for a   thin tectonic   plate  $\Omega$, see
 \cite{Reddy_2007}, and the  dynamics of  the plate described   by  a  perturbed  biharmonic wave  equation for vertical  displacement $u (x,t)$
 \begin{equation}
\label{perturbed_biharmonic}
H \rho u_{tt} + D \Delta^2 u + \nabla Q \nabla   u = 0,
\end{equation}
We formulate the  appropriate boundary conditions on $\partial \Omega \equiv \Gamma$ which are derived from  the  corresponding  Hamiltonian.
\medskip

In (\ref{perturbed_biharmonic}) the  flexural rigidity is 
\[ D = \frac{ H^3 E}{12(1-\sigma^2)} \equiv D_H  = D_1 \times  H^3  =  1.56 \times  H^3 \times 10^{10}\,\,  \frac{kg\, m^2}{sec^2}\]
 is  defined via 
 \begin{itemize}
 \item the Young modulus
 $E = 17.28 \times 10^{10} \,\,\frac{ kg }{m^{1}\, sec^{2}}$, 
 \item the Poisson  coefficient $\sigma = 0.28$, 
 \item the thickness $H \,\,\,\sim 3 \times 10^4 \,\,m - 10^5 \,\, m$,  and  
 \item the density $\rho =  3380 \frac{kg}{m^3}$ of the  plate.
 \end{itemize}

 The tangential tension   in the middle plane  of the plate  is modelled  by the  symmetric  elliptic  operator
\[ \nabla Q \nabla  u \equiv H T u \equiv \nabla Q_H \nabla  u = H \left[T_{x} u_{xx} +  2 T_{xy} u_{xy} + T_y u_{yy} \right]\equiv H \nabla Q_1 \nabla  u,\]
as in  \cite[Chapter 4]{Mikhlin_variational}, and constrained  in our case by the  maximal non-destructing  estimate  from above $  Q_1   \leq 0.3 \times 10^{10}$, see
\cite{FP_arxiv_2015}.

\medskip

 The boundary  conditions  for  the  biharmonic  wave equation (\ref{perturbed_biharmonic}) are derived  based  on the
  Hamiltonian
\begin{eqnarray*} 
 {\mathcal E} (u) & = &   \frac{1}{2}\int_{\omega}\left[H \rho u_t^2 +  D |\Delta u|^2 +  2 D(1-\sigma) [ |u_{xy}|^2 - u_{xx}u_{yy} ]+ \langle \nabla u, Q \nabla u\rangle \right]  d\Omega  \\ && \;\;\;\; +   \frac{1}{2} \int_{\Gamma}\beta \bigg|\frac{\partial u}{\partial n}\bigg|^2 d\Gamma, \label{Hamiltonian}
\end{eqnarray*}
 see our  Appendix 1,  and a more detailed discussion in  \cite[Chapters 4,8]{Mikhlin_variational}, as well as \cite{Birman_1956}, with regard  to the   boundary  bending defined  by   an  elastic
  bond  $\beta \frac{kg m} {sec^{2}}$.
  
  \medskip

  The  expression (\ref{Hamiltonian}) can be transformed, under the Dirichlet boundary condition  $u\bigg|_{\Gamma} = 0$ 
   into the following equation.
\begin{equation}
\label{Hamiltonian_Dirichlet}
 {\mathcal E}_D (u) =  \frac{1}{2}\int_{\omega}\left[H \rho u_t^2 +  D |\Delta u|^2 +
  \langle \nabla u, Q \nabla u\rangle \right]  d\Omega +  \frac{1}{2} \int_{\Gamma}[\beta - D \frac{1-\sigma}{r} ] \bigg|\frac{\partial u}{\partial n}\bigg|^2 d\Gamma,
\end{equation}
where $r$ is the curvature radius of the boundary (positive or negative depending on the position of the center of  the curvature),  see for instance  \cite{ Birman_1956}.

Minimizing  of the spacial  part  of  the  Hamiltonian (\ref{Hamiltonian_Dirichlet}) leads to the corresponding  ``natural''  boundary conditions:
\begin{equation}
 \label{bc_natural}
 [\beta - D \frac{1-\sigma}{r} ] \frac{\partial u}{\partial n} +  D \Delta u \bigg|_{\Gamma} = 0.
 \end{equation}
Mikhlin proved, see \cite{Mikhlin_variational}, that  the thin  Kirchhoff plate is  stable  if  $ [\beta - D \frac{1-\sigma}{r} ] \geq 0$   and $ Q \geq 0$,
  corresponding  to stretching  in the middle plane.  He also  considered  the  contracting  tension of  the  middle plane and  found  sufficient
conditions for  stability, see  \cite[Chapters 5,8]{Mikhlin_variational}.   Later Heisin,\cite{Heisin_1967}, actually noticed in an experiment that that  plates of
 ice  are  unstable  with respect to  certain  contracting  tension  in the middle plane.
In Section 3  we shall consider  an example  of  a  circular  plates  with   centrally  symmetric boundary conditions. The  corresponding  wave equations in
 the active zone $\Omega_{\varepsilon}$ are
  \begin{equation}
\label{biharmonic_varepsilon}
H_{\varepsilon}\rho w_{tt} +  D_{\varepsilon}  \Delta^2 w  +  Q_{\varepsilon} \Delta w  = 0,
\end{equation}
and on the complement $\Omega_c$.  Here the  tangent  contraction in the middle plane  is  neglected ($Q_{c} = 0 $), so that
  \begin{equation}
\label{biharmonic_c}
H_c\rho w_{tt} +  D_{c}  \Delta^2 w  +  Q_{c} \Delta w  = 0.
\end{equation}
With  properly selected  parameters  $H_{\varepsilon} = 3 \times 10^4\,\, m,\, H_c = 10^5 m,\, D_{\varepsilon} = D_1 \times H^3_{\varepsilon},\, D_c = D_1 \times H_c^3$\dots,
 solutions admit a spectral representation  constructed  with the use of Bessel functions.
 
 \medskip

 In next section we ensure that the resonance condition  $\lambda_{\varepsilon} = \lambda_c$ is satisfied
 for  a pair  of circular  tectonic plates. In  Section 4 we  impose a  weak bond  onto  a family of  oscillators with a  multiple  eigenvalue and  observe the dynamics   of the
 corresponding  perturbed  system.  It is this that  exposes  the  beat  phenomena involving the energy transfer between  oscillators   originally  constrained by the resonance condition.
 
Ultimately we will estimate the transition coefficient for the energy transfer for the simplest solvable model of  disjoint oscillators with relatively small masses and close
 frequencies $[\omega]$, perturbed  by  imposing  a  bond   of  them  with  an  oscillator  with  larger  mass  $M$  and  frequency  $\Omega \approx \omega$. 
 
This leaves the challenging problem of realising this  program for the  calculation  the energy transfer coefficient for oscillator's  to  the more interesting  system  of tectonic  plates, based  on fitted
 zero-range model.  We may discuss this elsewhere.
 
 \medskip

\section{Example : A  circular  active zone. }

Consider a thin  circular plate $\Omega$ divided  by a crack $\Gamma_{\varepsilon}$ into two   complementary  parts: the circular  active zone  $\Omega_{\varepsilon} $ and  the
 ring-like  complement $ \Omega_c$, centered  at  $\Omega_{\varepsilon} $. 
 
   \begin{figure}
[ht]
\begin{center}
\includegraphics [width=4in] {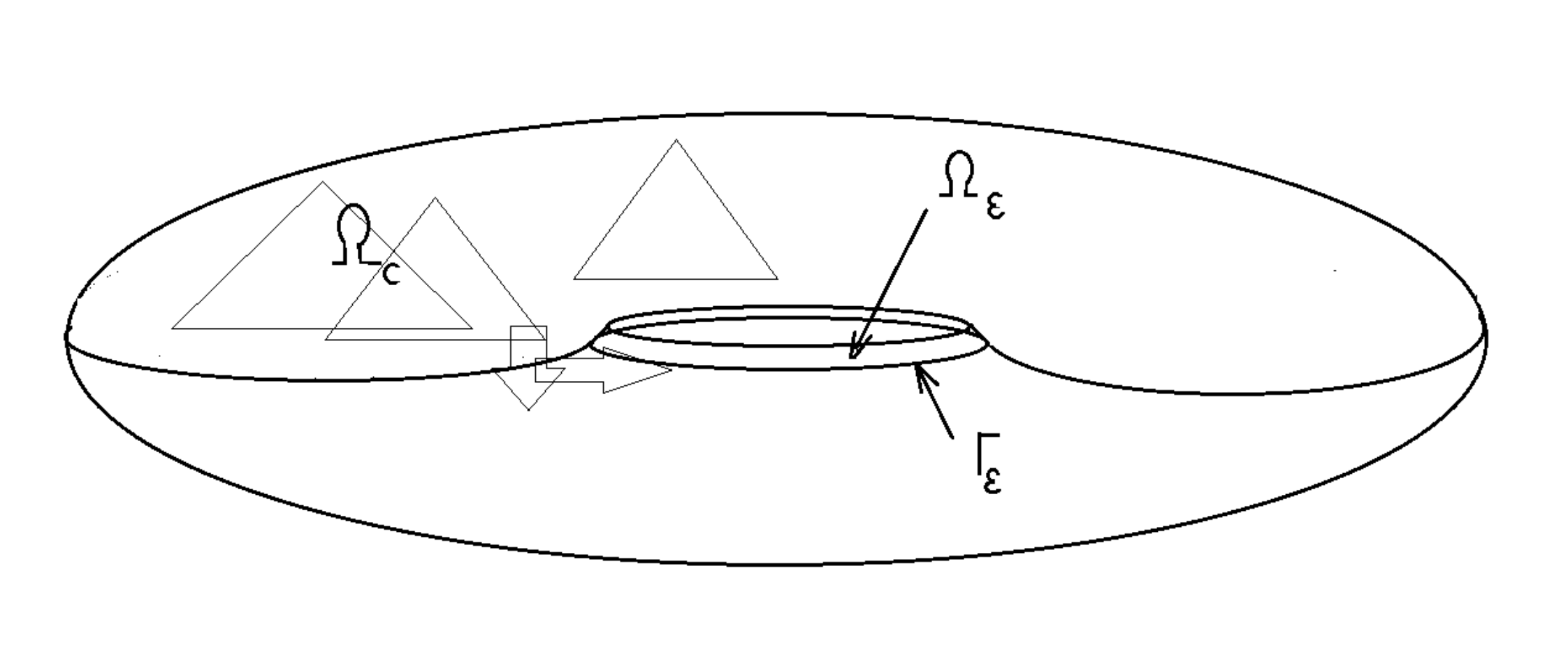}
\end{center}
\label{F:figure 4}
\caption
{The  small  circular  tectonic plate $\Omega_{\varepsilon}$ (the active zone) contacts  the complementary  large  circular tectonic plate $\Omega_c$  along
the boundary $\Gamma_{\varepsilon}$,
while  the  large plate  is loaded   and covers  the small  plate  on the  contact line $\Gamma_{\varepsilon}$.  Because of the load  and  the  special geometry
 of the contact,  the  large  plate develops a  normal stress (vertical  arrow) resulting in  bending  of the plate and a corresponding storage of  elastic energy. }
\end{figure}

\bigskip

We begin with the case  of  non-connected  parts, imposing  on both sides of $\Gamma_{\varepsilon}$ the
   kinematic boundary condition  $w_{\varepsilon}  \bigg|_{\Gamma} =  w_{c} \bigg|_{\Gamma} = 0 $ and independent  free-reclining boundary conditions, (see Appendix 1),   or  Neumann boundary conditions  imposed  on elements of the  domain of  the generators  $L_{\varepsilon}, L_{c} $ for the relevant  biharmonic  wave equations:

  \begin{equation}
\label{wave_biharmonic_varepsilon}
H \rho w_{tt} +  D_{\varepsilon}  \Delta^2 w  +  Q_{\varepsilon} \Delta w  = H \rho w_{tt} + L_{\varepsilon} w = 0,
\end{equation}
 with tangential  compression $ Q_{\varepsilon} \Delta \equiv T,\,\, Q_{\varepsilon} > 0$, characterized by the positive scalar $ Q_{\varepsilon}$ and an elastic
 bond $\beta$ applied
 on the boundary, as  in  (\ref{bc_natural}):
\begin{equation}
\label{bc_varepsilon}
w\bigg|_{\Gamma_{\varepsilon}} = 0,\,\,  \beta_{\varepsilon} \frac{\partial w}{\partial n} +   D_{\varepsilon}\Delta w  \bigg|_{\Gamma_{\varepsilon}}  -
\frac{ D_{\varepsilon} (1-\sigma)}{\varepsilon} \frac{\partial w}{\partial n} \bigg|_{\Gamma_{\varepsilon}} = 0.
\end{equation}

 On  the  complement   $\Omega_c$ we  neglect the  tangential compression, $Q_c = 0$, but do keep  the bending and the  elastic bond
 \begin{equation}
 \label{wave_biharmonic_c}
\rho H w_{tt} +  D_c  \Delta^2 w  = \rho H w_{tt} +  L_c   w  = 0,
\end{equation}
for the boundary condition  on the  outer side  of the  crack $\Gamma_{\varepsilon}$
\begin{equation}
\label{bc_c}
w\bigg|_{\Gamma_{\varepsilon}} = 0,\,\,  \beta_c \frac{\partial w}{\partial n} + D_c \Delta w  \bigg|_{\Gamma_{\varepsilon}}  -
\frac{D_c (1-\sigma)}{\varepsilon} \frac{\partial w}{\partial n} \bigg|_{\Gamma_{\varepsilon}} = 0,
\end{equation}
 and  for the  boundary condition on the remote part $\Gamma_a$ of the boundary  of ring-like complement  $\Omega_c$ we have
 \begin{equation}
\label{bc_a}
 \Psi_{c}\bigg|_{\Gamma_a} = 0,\,\,
\left[D_c\frac{1-\sigma}{r} - \beta_c\right]
\frac{\partial \Psi_{c}}{\partial n} - D_c \Delta \Psi_{c} \bigg|_{\Gamma_a} = 0.
\end{equation}
The  spectral  characteristics of   the  generators   $L_{\varepsilon}$ and $ L_c$ of the  wave dynamics on both $\Omega_{\varepsilon}, \Omega_c$, for given
geometrical and physical parameters can be  recovered by  separation of  variables.

\medskip

 We attempt  to select  the  parameters $D_{\varepsilon} =  H^3_{\varepsilon \times D_1}, \varepsilon,
 Q_{\varepsilon} = H_{\varepsilon} \times Q_1 $ with regard of  the resonance frequency $\nu_0 = 2 \times 10^{-4}$ Hz and  later choose the  outer radius $a$ of
  the complement $\Omega_c$ such
  that  one of  ground  frequencies  of the  biharmonic  generator $L_c$ also  coincides  with  $\nu_0 = 2 \times 10^{-4}$ Hz. Then  the  so constructed pair  of
  operators   $L_{\varepsilon}, L_c$
  can be perturbed by  the connecting  boundary  condition  on  the  crack, such that  the multiple eigenvalue $\omega^2_0 =  4\pi^2 \, \nu_0^2$  would split
  into a starlet
   $\omega^2_0 \longrightarrow \omega^2_0  \left(1 +  \delta [\alpha] \right),\,[\alpha] = [\alpha_1,\alpha_2]$ as  described below, implementing   the beating  of
    corresponding  spectral modes on $Q_{\varepsilon},Q_{c}$.

 All  data, except the radius  $\varepsilon$ of the active zone and  the outer radius  of the complement, are  selected as  in the  previous section,
 but we may  assume  the freedom of an adiabatic change, with time,  of the  tangent  tension in  the  middle plane
 as   $Q_{1} = q H\times 10^9,\,\, 1\leq q < 3$ below the destruction limit  of the  active zone $\Omega_{\varepsilon}$,
implying the change of  eigenfrequencies   $\omega_{\varepsilon} (q)$ of the active zone depending  on the  compressing tension.  

\medskip

 Removing  the common factor $H_{\varepsilon}$ from  the coefficients  of the  wave  equation  on the actice zone $Q_{\varepsilon}$   and from  the boundary  conditions  we
   obtain an equivalent  form of  (\ref{wave_biharmonic_varepsilon}) and the spectral problem on $\Omega_{\varepsilon}$ with  $H = H_{\varepsilon}$:

\begin{equation}
\label{spectral_varepsilon}
D_1 H^2 \Delta^2 w + Q_1 \Delta w = \omega^2 \rho w, \,\, \frac{\partial w}{\partial n}\bigg|_{\Gamma_{\varepsilon}} = u \bigg|_{\Gamma_{\varepsilon}} = 0.
\end{equation}
In our earlier work \cite{FP_arxiv_2015} we  estimated the  small  eigenvalues of the  above equations neglecting
the  boundary effects.  This  is acceptable  for  large  plates, see \cite{FP_arxiv_2015}, but  we notice that  improved
results  can be  obtained  based on a more accurate  spectral analysis  of  the biharmonic  generator  $L_{\varepsilon}$, using
a factorization of the above equation  with  Dirichlet-Neumann boundary conditions at  the  crack:

\begin{eqnarray}
0&=& \left[-\sqrt{D_1} H \Delta - \frac{Q_1} {2\sqrt{D_1} H } - \sqrt{\Omega^2 \rho  +  \frac{ Q_1^2}{4 D_1 H^2}}\right] \nonumber \\ && \times
\left[-\sqrt{D_1} H \Delta - \frac{Q_1} {2\sqrt{D_1} H } + \sqrt{\Omega^2 \rho  +  \frac{ Q_1^2}{4 D_1 H^2}}\right] \label{factorization}
\end{eqnarray}

For the circular  active zone  $\Omega_{\varepsilon}$  and  centrally symmetric  $ w = \Psi_{\varepsilon} $,\,\,  $\Delta \equiv \Delta_0$ is the corresponding  radial Laplacian,
hence the  above equation (\ref{factorization}) has a solution vanishing on  $\Gamma_{\varepsilon }: r = \varepsilon$  presented as  a linear combinations  of  Bessel functions
and  modified Bessel functions  with  appropriate arguments:
\begin{equation}
\label{bessel_solutions}
\Psi_{\varepsilon} (r) =  \frac{J_0 \left(\left[ \frac{\omega^2 \rho }{H^2 D_1} e^{2\Theta}  \right]^{1/4} r \right)}
{J_0 \left(\left[ \frac{\omega^2 \rho }{H^2  D_1} e^{2\Theta}  \right]^{1/4} \varepsilon \right)} -
\frac{I_0 \left(\left[ \frac{\omega^2 \rho }{H^2 D_1} e^{-2\Theta}  \right]^{1/4} r \right)}
{I_0 \left(\left[ \frac{\omega^2 \rho }{H^2 D_1} e^{-2\Theta}  \right]^{1/4} \varepsilon \right)}.
\end{equation}

Here  $\sinh \Theta = \frac{Q_1}{2\omega H\sqrt{D_1 \rho}} $ reveals  the dependence  of the resonance on the tangential tension  $Q$.  Besides
\begin{eqnarray*}
-\Delta J_0 \left(\left[ \frac{\omega^2 \rho}{D} e^{2\Theta}  \right]^{1/4} r \right) & = &   \left[ \frac{\omega^2 \rho}{D_1 H^2} e^{2\Theta}  \right]^{1/2}\,
J_0 \left(\left[ \frac{\omega^2 \rho}{D} e^{2\Theta}  \right]^{1/4} r \right),\\
\Delta I_0 \left(\left[ \frac{\omega^2 \rho }{D} e^{-2\Theta}  \right]^{1/4} r \right) & = &   \left[ \frac{\omega^2 \rho }{D_1  H} e^{-2\Theta}  \right]^{1/2}\,
I_0 \left(\left[ \frac{\omega^2 \rho }{D_1 H^2} e^{-2\Theta}  \right]^{1/4} r \right)
\end{eqnarray*}

We are interested in  small positive  eigenvalues $\omega^2$ of the equation (\ref{spectral_varepsilon}) which may be in resonance with the lower
eigenmodes of  the  complementary part $\Omega_c$ of the plate.  To  recover an algebraic equation for the eigenfrequencies we  substitute the spectral parameter  $\omega$ by
a new  spectral parameter $\Theta$ connected to $\omega$ by the equation 
\[ \sinh \Theta  = \frac{Q}{2\omega \sqrt{D  \rho}}  \]
or 
\[ \omega  \equiv \frac{Q}{2  \sinh \Theta \sqrt{D  \rho}} = \frac{Q}{ \left[e^{\Theta} - e^{-\Theta}\right]\sqrt{D  \rho}}, \Theta > 0.\]
 The unperturbed spectral  problem
with the new spectral parameter  is defined by the Dirichlet-Neumann  boundary condition  while constructed of Bessel functions as

\begin{equation}
\label{Dir_sol_varepsilon}
\Psi_{\varepsilon} (r) =  \frac{J_0 \left(e^{\Theta/2}\frac{\sqrt{Q}}{\sqrt{D}\sqrt{e^{\Theta} - e^{-\Theta}}} r \right)}
{J_0 \left(e^{\Theta/2}\frac{\sqrt{Q}}{\sqrt{D}\sqrt{e^{\Theta} - e^{-\Theta}}} \varepsilon \right)} -
\frac{I_0 \left(e^{-\Theta/2}\frac{\sqrt{Q}}{\sqrt{D}\sqrt{e^{\Theta} - e^{-\Theta}}} r \right)}
{I_0 \left(e^{-\Theta/2}\frac{\sqrt{Q}}{\sqrt{D}\sqrt{e^{\Theta} - e^{-\Theta}}} \varepsilon \right)},
\end{equation}
 and   satisfies   zero boundary condition identically   on the inner side of the crack $\Gamma_{\varepsilon}$.
The Neumann boundary condition on the inner side of the crack 

\begin{eqnarray} \nonumber
 \frac{\partial \Psi_{\varepsilon}}{\partial n} (\varepsilon) & = & 
 e^{\Theta/2}\frac{\sqrt{Q_1}}{\sqrt{D_1 H^2}\sqrt{e^{\Theta} - e^{-\Theta}}}
 \frac{J'_0 \left(e^{\Theta/2}\frac{\sqrt{Q_1}}{\sqrt{D_1 H^2}\sqrt{e^{\Theta} - e^{-\Theta}}} \varepsilon \right)}
{J_0 \left(e^{\Theta/2}\frac{\sqrt{Q_1}}{\sqrt{D_1 H^2}\sqrt{e^{\Theta} - e^{-\Theta}}} \varepsilon \right)} \\&& -
 e^{-\Theta/2}\frac{\sqrt{Q_1}}{\sqrt{D_1 H^2}\sqrt{e^{\Theta} - e^{-\Theta}}}
\frac{I_0 \left(e^{-\Theta/2}\frac{\sqrt{Q_1}}{\sqrt{D_1 H^2}\sqrt{e^{\Theta} - e^{-\Theta}}} \varepsilon \right)}
{I_0 \left(e^{-\Theta/2}\frac{\sqrt{Q_1}}{\sqrt{D_1 H^2}\sqrt{e^{\Theta} - e^{-\Theta}}} \varepsilon \right)} = 0\label{Neumann_inner_0}
\end{eqnarray}
and defines  the spectrum of the  unperturbed   Dirichlet-Neumann  problem  for selected values  of the  parameters  involved.

 We  select the  typical  geometrical and physical
  parameters of the active zone $\Omega_{\varepsilon}$ as 
  \begin{itemize}
  \item $ D_1 \times H^2 = 1.56 \times 10^{10} \times H^2 \frac{kg m}{sec^2} $, 
   \item $  H = H_{\varepsilon} \sim 3 \times 10^4 m $, 
   \item $\rho = 3380 \,\frac{kg}{m^3}$,  $ \sigma = 0.28$, $Q \approx 3\times 10^9 \frac{kg }{m \,\sec}$, and
   \item $\varepsilon =  {\rm radius} \Omega_{\varepsilon}\sim 2.6 \times 10^5 m$
   \end{itemize}
   (with a  little  bit of  accurate tuning).  Then
     we are able to substitute  the Bessel functions $J_0$ and the corresponding derivatives by the  asymptotics  for ``large'' values of the argument,
      and  $I_0$ and the corresponding derivatives  by the asymptotics  for small ones, with regard  of  the factor  $e^{-\Theta}$, calculated  for
       the  selected   resonance  frequency  $ \nu = 2 \times 10^{-4}$  Hz of the lower  eigenmodes  of large  tectonic plates.  
\[
\sinh \Theta = \frac{Q}{4\pi \nu \sqrt{\rho D}} =  \frac{3 \times 10^9}{12.56 \times 2 \times 10^{-4} \times
 3 \times 10^4 \times 1.26 \times 10^5}=  5.5,\,\mbox{and}\,\, e^{\Theta} = 11.
\]
Then the arguments of  the Bessel functions and  their derivatives  can be calculated  to be
\begin{eqnarray}\nonumber
J_0 (\omega^{1/2} \left(\frac{\rho}{D})^{1/4}\, e^{\Theta/2} \varepsilon \right) & = & 
 J_0 \left(\left[\frac{12.6 \times 10^{-4}\times 58\times 11 \times 6.76 \times 10^{10}}
 {1.26 \times 10^5 \times 3\times 10^4} \right]^{1/2}\right)\nonumber \\ & \approx &J_0 (\frac{\sqrt{1505}}{10}) = J_0 (3.9), \\
 I'_0 (\omega^{1/2} \left(\frac{\rho}{D})^{1/4}\, e^{\Theta/2} \varepsilon \right)& = & 
\label{arguments}
 I_0 \left(\left[\frac{12.6 \times 10^{-4}\times 58 \times 6.76 \times 10^{10}}
 {1.26 \times 10^5\times 11 \times 3\times 10^4} \right]^{1/2}\right)\nonumber \\ & \approx & I_0 (\frac{\sqrt{61}}{10}) \approx I_0 (0.8)
\end{eqnarray}

Hereafter we use  standard Taylor  asymptotics  of the modified Bessel function  $I_0$ for small values of  argument ($< 1 $) and the exponential asymptotics
 for ``large" arguments ($ \geq 3.9 $) in  $J_0$,
\[
J_0 (z) \approx \frac{\cos (z - \pi/4) + O (1/z)}{\sqrt{\frac{\pi z}{2}}},\,\,
J'_0 (z) \approx \frac{-\sin (z - \pi/4) + O (1/z)}{\sqrt{\frac{\pi z}{2}}}
\]
\begin{equation}
I_0 (z) \approx \frac{\cosh z  + O (1/z)}{\sqrt{\frac{\pi z}{2}}},\,\,
I'_0 (z) \approx \frac{\sinh z  + O (1/z)}{\sqrt{\frac{\pi z}{2}},
}.
\end{equation}

Then  we  notice that the equation  (\ref{Neumann_inner_0}) is  satisfied, with data, selected above, up to
an error $\sim 0.1$. This means that  our guess concerning  the  magnitude  of  the frequency was reasonably accurate  for an  active zone  with radius  $2.6 \times 10^5$
 with  standard  physical characteristics and circular shape.

A more  profound  correspondence between the  interval  $150 \mu Hz - 250 \mu Hz$ of typical  frequencies and  various shapes of the  active zone requires a
further analysis of  the corresponding  dispersion equation (analog of \ref{Neumann_inner_0}) and is postponed.

\bigskip

\section{Resonance conditions for  circular plates.}

To reveal the resonance condition for  the circular  plate divided  into two parts: the circular active zone $\Omega_{\varepsilon}$   and the ring-like
complement  $\Omega_c  = \left\{ \varepsilon <  r < a \right\}$, which is
centered at  $\Omega_{\varepsilon}$ and  elastically  disconnected
from the  active zone due to  independent Dirichlet-Neumann conditions on both sides  of the  common boundary
$\Gamma_{\varepsilon}$ we must select the geometrical  parameters  of the  complement  such that the
 disconnected  spectral  problem has  a multiple eigenvalue  $\omega^2_0 = 4\pi^2 \nu_062$.  Then  the
 perturbation of the  disconnected spectral  problem defined by replacing the disconnecting  boundary  conditions  by  an interactive condition  would
  reveal  a  splitting  of the multiple eigenvalue, and, eventually, the the resonance  beating of SGO  modes, localized  on  $\Omega_{\varepsilon}, \Omega_c$.  This would
   imply  migration of energy between the locations.
 It is  technically  convenient to  consider a circular  pate $\Omega'_c$  with the same outer radius  $a$, but without a hole  reserved  for $\Omega_{\varepsilon}$ at
  the center.  The  perturbation obtained by replacement  $\Omega_c \to \Omega'_c$ should not affect the  part of spectrum  corresponding  to the standing waves  in which
  lengths  exceed  the  geometric size  of the details
 affected  by the  change, in our case  the  radius $2.6 \times 10^5\,\, m $ of the  hole $\Omega_{\varepsilon}$ (the active zone). This condition
 is obviously  satisfied  for ground  SGO  modes  on  the  disc $\Omega'_c$ with $ R  \approx  5000 m$.  Now estimation of  the eigenfrequencies  of  the  complement  for an  equivalent  circular  plate  $\Omega'_c$, with  radius   $a \approx 5 \times 10^6 m $,
 neglecting  relatively  small terms, compared  with the typical ground flexural  wavelengths $2a \approx 5\times 10^6 m $, the  hole radius  $\varepsilon$. We  assume  for  $\Omega_c,\Omega'_c : D_a = 1.56 \times 10^{10}\times H_a^2 \frac{kg m}{sec^2} $, $  H_a = H_{c} \sim 10^5 m $,
 $Q_c \approx 0 $ and $\rho = 3380 \,\frac{kg}{m^3}$.
We also use  of the asymptotics  of the derivatives  of the Bessel  functions for the  solutions  $\Psi_a$ of the spectral problem
\begin{equation}
\label{spectral_a}
D_a \Delta^2 w  = \omega^2 \rho w
\end{equation}
on the complement, with $Q_c = 0,\, \Theta_c = 0$:

 \begin{equation}
\label{bessel_solutions_a}
\Psi_{a} (r) =  \frac{J_0 \left(\left[ \frac{\omega^2 \rho }{D_1 H^2_a} \right]^{1/4} r \right)}
{J_0 \left(\left[ \frac{\omega^2 \rho }{D_1 H_a^2}   \right]^{1/4} a \right)} -
\frac{I_0 \left(\left[ \frac{\omega^2 \rho }{D_1 H^2_a}  \right]^{1/4} r \right) }
{I_0 \left(\left[ \frac{\omega^2 \rho }{D_1 H^2_a}  \right]^{1/4} a \right)},
\end{equation}

For the  Neumann dispersion equation on the remote part  $\Gamma_{a}$  of the boundary
\begin{equation}
\label{dispers_neumann}
\frac{\partial \Psi_{a}}{\partial n} (\varepsilon) \times {\sqrt{\frac{\pi z}{2}}}
 =
 \frac{J'_0 \left(\sqrt{\omega} \left(\frac{\rho}{D_1 H^2_a} \right) a \right)}
{J_0 \left(\sqrt{\omega} \left(\frac{\rho}{D_1 H^2_a} \right) a \right)} -
\frac{I'_0 \left(\sqrt{\omega} \left(\frac{\rho}{D_1 H^2_a} \right) a \right)}
{I_0 \left(\sqrt{\omega} \left(\frac{\rho}{D_1 H^2_a} \right) a \right)} = 0.
\end{equation}

Using  the asymptotics  of  Bessel functions for large  arguments  on the remote part of the boundary $r=a$ we find that
\begin{eqnarray*}
\frac{I'_0 \left(\sqrt{\omega} \left(\frac{\rho}{D_1 H^2_a} \right) a \right)}{I_0 \left(\sqrt{\omega}
 \left(\frac{\rho}{D_1 H^2_a} \right) a \right)}\approx 1 &=& - \tan  \pi (l  -\frac{1}{4}), \;\;\; 
\frac{J'_0 \left(\sqrt{\omega} \left(\frac{\rho}{D_1 H^2_a} \right)^{1/4}  a \right)} {J_0 \left(\sqrt{\omega} \left(\frac{\rho}{D_1 H^2_a} \right)^{1/4} a \right)}
 \\ & \approx& - \tan
\left(0.6 a \times 10^{-7} - \frac{\pi}{4} \right).
\end{eqnarray*}

 This results  in $0.6  a \times 10^{-7}- \pi/4 = \pi l  - \pi/4$ and hence $a \approx 5 \times 10^6 m$, for $l=1$
 in agreement with our preliminary guess.

 More detailed estimation of the  radius of the  ring - like complement $\Omega_c \equiv \left(\varepsilon < r < a \right)$, centered at the active zone
 $\Omega_{\varepsilon}$ may be derived from an  explicit construction of the basic Bessel solutions and  the  Neumann-to-Dirichlet  map  of the  biharmonic
 d'Alambert equation (\ref{spectral_a}) on the ring with  and appropriate symmetry, compatible with typical  standing waves  on the complement.
\begin{eqnarray}\nonumber
J_p (z) &\approx& \frac{\cos (z - p\pi/2 - \pi/4)}{\sqrt{\pi z/2}} + O(|1/z|) e^{|\Im z|}, \\
J'_p (z) &\approx& \frac{-\sin (z - p\pi/2 - \pi/4)}{\sqrt{\pi z/2}} + O (|1/z|) e^{|\Im z|} \nonumber \\
I_p (z) &\approx& \frac{e^z}{\sqrt{\pi z/2}} + O (e^{|{\Im z}|}),\\
I'_p (z) &\approx & \frac{e^z}{\sqrt{\pi z/2}} + O(e^{|{\Im z}|}) 
\end{eqnarray}

One can also consider  a circular  active zone  on  an arbitrary  plate  with smooth boundary, see  \cite{Timoshenko_circular}, or derive  the formulae  for
  the  ND -map  based  on  integral equation techniques,
see  \cite{Birman_1956}.

Henceforth we wish to consider sectorial  circular  active zones which admit a separation of variable along with  solutions  of the relevant  perturbed  biharmonic equation
\begin{equation}
\label{biharmonic_p}
D \Delta^2 u  +  Q \Delta u = \omega^2 \rho u
\end{equation}
 combined on  the sector $\Omega^p_{\varepsilon} :  0\leq r < \varepsilon, \,\, 0 < \varphi < \pi / p $ of above  Bessel  functions with index $p$  and
  the  corresponding circular harmonics  $\cos p \phi,\,\sin p \phi  $  on the sector, as above.

  The  roles of  basic regular  solutions of the biharmonic d'Alambert equation play the products  of  circular harmonics $\sin p \phi,  \cos p \phi $, with
  relevant  Bessel functions  ${\mathcal J}_{p} (r\sqrt{\omega} \left(\frac{\rho}{D} \right)^{-1/4} e^{\Theta/2})$,
  and the modified Bessel functions 
  \[ {\mathcal I}_{p} (r\sqrt{\omega} \left(\frac{\rho}{D} \right)^{-1/4} e^{-\Theta/2}).\]
 The  parameter $\Theta$ is  derived from the  factorization  of  the  d'Alambert equation as at (\ref{factorization}). Then, for $p \geq 1 $
 the  only continuous at $r=0$ solution are square integrable.  Hence continuous  solutions of the d'Alambert  equation, vanishing  on the circular
  part $\Gamma_{\varepsilon}$ of the  boundary can b  obtained as  linear combinations.
  \begin{equation}
  \label{psi_sin}
  \Psi_{s} (r) = \sin p \phi \left[\frac{J_{p} \left(\left[ \frac{\omega^2 \rho }{D} e^{2\Theta}  \right]^{1/4} r \right)}
{J_{p} \left(\left[ \frac{\omega^2 \rho }{D} e^{2\Theta}  \right]^{1/4} \varepsilon \right)} -
\frac{I_{p} \left(\left[ \frac{\omega^2 \rho }{D} e^{-2\Theta}  \right]^{1/4} r \right)}
{I_{p} \left(\left[ \frac{\omega^2 \rho }{D} e^{-2\Theta}  \right]^{1/4} \varepsilon \right)}\right],
  \end{equation}

   \begin{equation}
  \label{psi_cos}
  \Psi_{c} (r) = \cos p \phi \left[\frac{J_{\gamma} \left(\left[ \frac{\omega^2 \rho }{D} e^{2\Theta}  \right]^{1/4} r \right)}
{J_{p} \left(\left[ \frac{\omega^2 \rho }{D} e^{2\Theta}  \right]^{1/4} \varepsilon \right)} -
\frac{I_{p} \left(\left[ \frac{\omega^2 \rho }{D} e^{-2\Theta}  \right]^{1/4} r \right)}
{I_{p} \left(\left[ \frac{\omega^2 \rho }{D} e^{-2\Theta}  \right]^{1/4} \varepsilon \right)}\right].
  \end{equation}
  They satisfy   the Dirichlet boundary conditions $u\bigg|_{\Gamma} =
\Delta u\bigg|_{\Gamma} = 0 $ for $ \Psi_{s}$ and  the Neumann boundary condition $\frac{\partial u}{\partial n}\bigg|_{\Gamma} =
\frac{\partial \Delta u}{\partial n}\bigg|_{\Gamma} = 0 $ respectively.

\medskip

The spectral problem with Dirichlet-Neumann boundary condition
 $u \bigg|_{\Gamma_{\varepsilon}} = \frac{\partial u}{\partial n} \bigg|_{\Gamma_{\varepsilon}}$ on the circular part of the boundary is given by  linear combinations ({\ref{psi_sin},{\ref{psi_cos}), and  which satisfy  the corresponding  dispersion equations.
\begin{eqnarray}\nonumber
0& = & \frac{\partial \Psi_{s}}{\partial n} (\varepsilon)\, \sin \gamma \phi \left[ e^{\Theta/2}\frac{\sqrt{Q}}{\sqrt{D}\sqrt{e^{\Theta} - e^{-\Theta}}}
 \frac{J'_{\gamma} \left(e^{\Theta/2}\frac{\sqrt{Q}}{\sqrt{D}\sqrt{e^{\Theta} - e^{-\Theta}}} \varepsilon \right)}
{J_{\gamma} \left(e^{\Theta/2}\frac{\sqrt{Q}}{\sqrt{D}\sqrt{e^{\Theta} - e^{-\Theta}}} \varepsilon \right)} \right. \\
&& \left.- e^{-\Theta/2}\frac{\sqrt{Q}}{\sqrt{D}\sqrt{e^{\Theta} - e^{-\Theta}}}
\frac{I'_{\gamma} \left(e^{-\Theta/2}\frac{\sqrt{Q}}{\sqrt{D}\sqrt{e^{\Theta} - e^{-\Theta}}} \varepsilon \right)}
{I_{\gamma} \left(e^{-\Theta/2}\frac{\sqrt{Q}}{\sqrt{D}\sqrt{e^{\Theta} - e^{-\Theta}}} \varepsilon \right)}\right]. \label{Neumann_inner}
\end{eqnarray}
or a similar equation for $\Psi_{c}$. 

For sectors   characterized by $p \in (0,1)$ there are singular square-integrable solutions  $J_{-p} (z),\,\, I_{-p} (z)$ of the biharmonic  d'Alambert  equation,
see \cite{Sectorial_defect}. Then  we are able to consider the  singular (discontinuous) square integrable solutions of the biharmonic d'Alambert equation,
 considering them as  elements of the corresponding defect, see  \cite{Sectorial_defect}, and  construct self-adjoint extensions of the corresponding
 biharmonic operator   in the space of square-integrable functions  in an  sectorial  active zone on the boundary. The  corresponding operator extension
  machinery can be developed  with use  of extension procedure, similar to  one developed  above for the inner  zero-range active zone.
  
  \medskip

 When varying the tension parameter $Q$,  we  will come to the moment when  the  tangent  compression
  in the middle plane  is large enough  for  the  minimal  eigenvalue  on the active zone  of an unperturbed problem to
coincide  with an  eigenvalue  of the biharmonic  spectral problem ($Q = 0$) on  the complement.
At this point substituting the asymptotics of the Bessel functions $J_p, I_p$ for large and small values of the arguments, one can  obtain an  estimation for
 the contracting tension  $Q$,  this leads to the resonance conditions  for the  operators  $L_{\varepsilon}, L_c$  on  $\Omega_{\varepsilon}, \Omega_c$,
 without an interaction between  them.  We now explore this in the simplest case.

\medskip

\section{ A simple   model  of alternation.} 

Now we consider  a  simple universal  interaction  depending  on a small parameter between  similar  operators
  constructed  for the zero-range active zone. Based on this construction, we  observe the  relevant beating phenomenon and  estimate the  transfer
 coefficient. 
 
 \medskip 

The  system of  two  tectonic plates considered in the previous section is a special case  of a  decoupled oscillator's system  under
 the resonance conditions.   The  interaction between the  plates could be introduced by imposing free reclining  or natural
 boundary conditions.  
 
 We now make some simplifying assumptions:   we consider a  weakly coupled oscillatory system
 obtained  by  attaching  one (supposedly ``large'') multi-dimensional  oscillator $ X = (X_1, X_2,X_3, \dots X_{\mu}) \in  C_{\mu}$
 to  a  1D (supposedly ``small'')  oscillator characterized  by the coordinate  $ x = x_1 \in  C_{1}$.  
 
 The dynamics of the  resulting
  system is  defined by the system of  linear equations
\begin{equation}
\label{coupled_1}
\begin{array}{c}
m x_{tt} + v x +  b^+ X = 0,\\
M X_{tt} +  V X +  b x = 0.
\end{array}
\end{equation}
Here  $ m = m_1, v = v_1, M = [M_1, M_2, M_3\dots M_{\mu}], V = [V_1, V_2, V_3\dots V_{\mu}]$,  are  positive diagonal matrices acting in  the Hilbert spaces $ C_{1},\, C_{\mu}$,  $ C_{1} \stackrel{b} \longrightarrow C_{\mu},\,  C_{\mu} \stackrel{b^+} \longrightarrow C_{1} $.
 Interaction of the  oscillators  is  introduced  by  the Hermitian matrix $B : C_{1} \oplus C_{\mu} \longrightarrow C_{1} \oplus C_{\mu}\equiv K$
\begin{equation}
\label{interaction}
B =
 \left(
\begin{array}{cc}
0 & b^+\\
b & 0
\end{array}
\right) \equiv {\rm antidiag} \left(b^+, b \right),
\end{equation}
which plays  the role of  bonds imposed on the  boundary  data  of   solutions of the biharmonic equation  on the  border  of  the tectonic plate.  Separating  time, we  obtain the spectral  problem corresponding to the above  wave equation (\ref{coupled_1})
\[
v x   + b^+ X =  m\lambda x
 \]
 \begin{equation}
 \label{spectral_1}
 V X +  b x =  M  \lambda X.
 \end{equation}
 To  make this consistent  with the above model  of  tectonic plates  supplied with  natural  structure of active zones  we  assume 
   that  the  unperturbed
 multidimensional  oscillator  has a family of  unperturbed  eigenvalues   (square frequencies)
  $[\omega^0]^2 = v/m, \left[\Omega^0_s\right]^2 = V_s / M_s, s= 1,2,\dots \mu $. For  non-zero
interaction $ B \neq 0$ the eigenvalues   $\lambda^b_r =  (\omega_r^b)^2$
(and   the corresponding eigenfrequencies)  of  the  perturbed
selfadjoint  spectral  problem  are  found  from  the algebraic equation
\begin{equation}
\label{spectrum}
\left [ \lambda m - v  \right] a + b^+ \frac{I}{ V - M \lambda} b  a \equiv {\mathcal{M}} (\lambda) a  = 0,
\end{equation}
obtained via the elimination of $X$  from the second equation  in  (\ref{spectral_1}).  In particular,  for the 1D case,  $\mu = 1$, we have  two  unperturbed  frequencies  
\[ v/m = \omega^0 \equiv (\omega^0_1)^2,\, V/M = \Omega^0 \equiv (\omega^0_2)^2\]
and a quadratic equation for the perturbed  eigenvalues 
\[ \lambda^b_r = (\omega^b_r)^2,\, r = 1,2, \]
while the unperturbed are $\lambda^0_1 = (\omega^0)^2,\,
 \lambda^0_2 = (\Omega^0)^2 $.  Denoting  $\bar{\lambda}^0 = \frac{(\omega^0)^2 + (\Omega^0)^2}{2},\,\,\delta \lambda^0 = \bigg| \frac{(\omega^0)^2 - (\Omega^0)^2}{2}\bigg| $, we  find the perturbed  frequencies/eigenvlues  $\lambda^b_r \equiv (\omega^b_r)^2$ from the quadratic equation
\[
\lambda^2 -  2  \lambda \,\delta\lambda^0 +  \lambda^0_1 \,\lambda^0_2 = \frac{b^2}{m M}. \]
This gives
\[ \lambda^b_r =  \bar{\lambda}^0  \pm  h, \,\, \mbox{with} \,\, h^2 = (\delta\lambda)^2 + \frac{b^2}{mM},\]
with $m_1 = m,\,\, m_2 = M$.
We further  assume that  $\frac{b^2}{m_1 m_2} << \delta^2 \lambda $. This assumption  allows us to  calculate
the perturbed  eigenvalues  $\lambda^b_{1,2}$  with  $h = \delta\lambda +  \frac{b^2}{2 m_1\, m_2\, \delta \lambda} $  approximately.  These turn out to be
\begin{eqnarray}
 \nonumber
\lambda^b_{1,2} &=&  \lambda^0_{1,2} \pm \frac{b^2}{2 m_1\, m_2 \, \delta \lambda},\,\,
 \omega^b_{1,2} = \omega^0_{1,2} \pm \frac{b^2}{4 m_1\, m_2\, \delta \lambda\, \omega^0_{1,2}} \\
 & = &  \omega^0_{1,2}\left[1 \pm \frac{b^2}{4 m_1\, m_2\, \delta \lambda^2 }\right]\frac{\delta \lambda}{\lambda_{1,2}^0}\approx
 \omega^0_{1,2} \frac{\delta\lambda}{\lambda_{1,2}^0}.\label{eigenvalues_perturbed}
\end{eqnarray}

Generally, for $\mu \geq 1 $, the  eigenvalues of the ultimate spectral  problem
$\lambda =  \lambda^b_1, \, \lambda^b_2, \dots \lambda^b_s,\dots \lambda^b_{1 + \mu}$ may be found  from the corresponding  determinant condition.  The  components   of the corresponding   eigenvectors $ \left\{a,\Psi^b_s\right\} \equiv
{\bf \Psi}^b_s$  in $K, C_{\mu}$  are  calculated  from the equation
\begin{equation}
\label{eigenvectors}
 \Psi^b_s = \left(V -  M \lambda^b_s \right)^{-1} b a,
\end{equation}
with  an appropriate normalization  $m|a|^2 \left[ 1 + \langle M \frac{I}{ V - M \lambda} b, \frac{I}{ V - M \lambda} b\rangle  \right] = 1$.
 If  the  interaction  $b$ is real and  the initial data of the  relevant Cauchy problem,   see below (\ref{Cauchy}), are real, then the corresponding  solution is real too. Similarly the  solution of the corresponding  inhomogeneous  equation,  with  a  real  function in the  right side,  is real as well.

\begin{figure}
[ht]
\label{Spectrum perturbed}
\begin{center}
 \includegraphics [width=5 in]{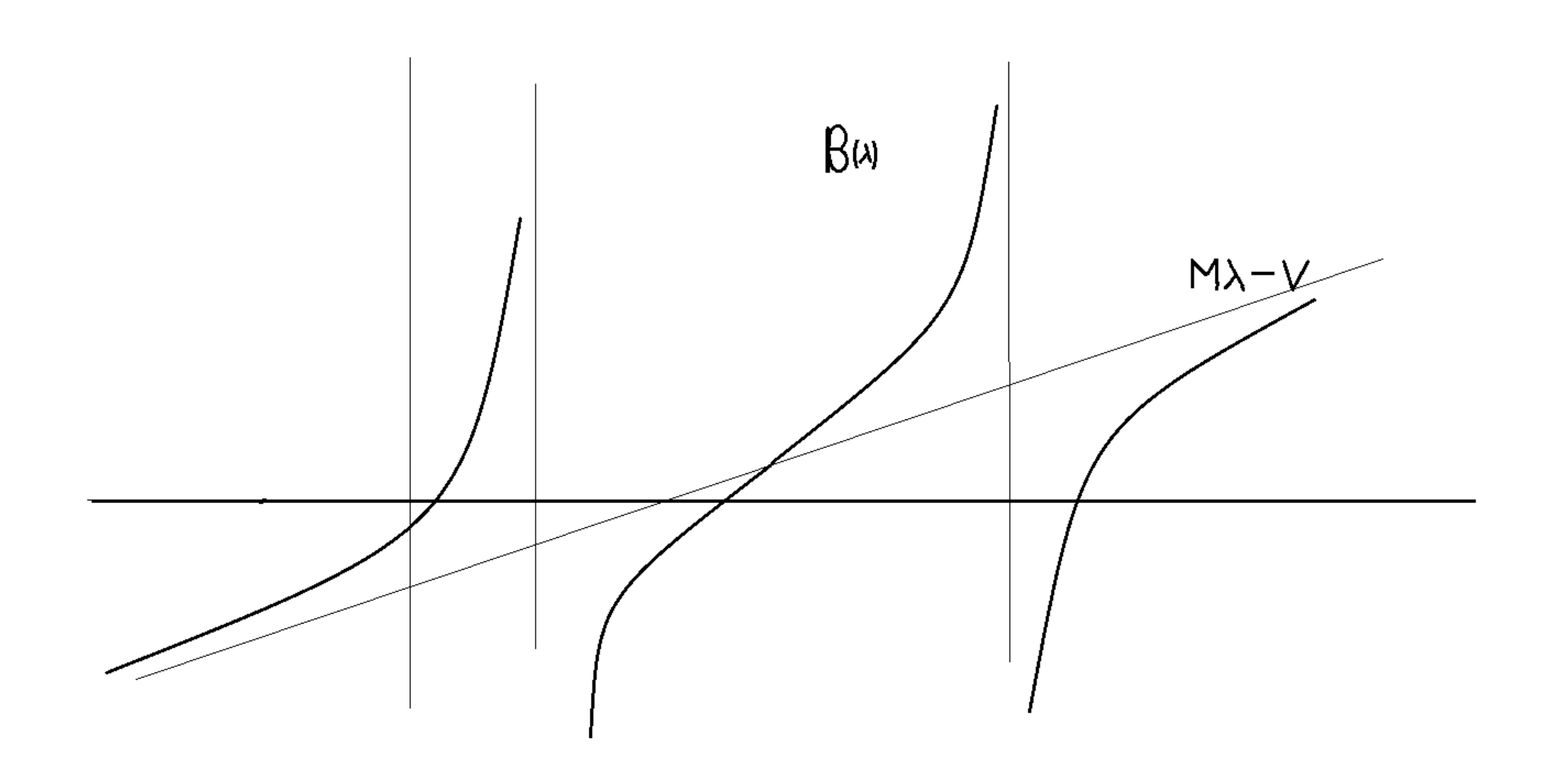}
\end{center}
\caption{In the more general case where  $K = C_1 \oplus C_{\mu},\, \mu \geq 1 $, the eigenvalues  $\lambda^b_s = (\omega^b_s)^2 $ are  found  as  graphical solutions based  on
 the  diagram, corresponding to  the algebraic equation  (\ref{spectrum}),  see \cite{KP_1974}.}
\end{figure}

In the simplest situation where $\mu = 1, K = C_1 \oplus C_1$ the  equation (\ref{spectral_1}) for $(u, U) = u_1, u_2$ and the unperturbed  frequencies
$ \omega^0 \equiv \omega^0_1, \Omega^0 = \omega^0_2 $  and  
\[ \bar{\lambda} = \frac{\lambda^0_1 + \lambda^0_2}{2},\,\, \delta \lambda = \frac{\lambda^0_1 - \lambda^0_2}{2}\]
 the equation (\ref{spectral_1}) can be represented as
\begin{equation}
\label{spectral_11}
\left(
\begin{array}{cc}
m_1\lambda^0_1 &  b^+ \\
b & m_2 \lambda^0_2
\end{array}
\right)
\left(
\begin{array}{c}
u_1\\
u_2
\end{array}
\right) = \lambda
\left(
\begin{array}{cc}
m_1 & 0 \\
0& m_2
\end{array}
\right)
\left(
\begin{array}{c}
u_1\\
u_2
\end{array}
\right),
\end{equation}
and  the  eigenvectors
\begin{equation}
{\bf{\Psi}}^b_1 = \left(
\begin{array}{c}
a_1\\
\frac{a_1b m_2^{-1}}{-\delta \lambda  - h}
\end{array}
\right),\,\,
 {\bf{\Psi}}^b_2 = \left(
\begin{array}{c}
a_2\\
\frac{a_2 b m_2^{-1}}{-\delta \lambda  + h}
\end{array}
 \right)
 \end{equation}
are  orthogonal  with respect  to $diag\,\,\,(m_1,m_2)$ and  normalized in $l_2 (m_1,m_2)$ with  $a_{1,2} =
\sqrt{\frac{h \pm \delta \lambda}{2h}}$. When discussing  the  general case $\mu \geq 1$
we may assume that the  eigenfunctions  ${\bf \Psi}^b_s = \left(\psi^b_s,\, \Psi^b_s \right)$
of the perturbed spectral  problem  (\ref{spectral_1}) correspond to simple eigenvalues and  are
orthogonal and normalized. Using  the eigenvectors,  the normal modes of the wave equation
\begin{equation}
\label{Cauchy}
\begin{array}{c}
m u_{tt} + v u  +  b^+ U = 0,\\
M U_{tt} +  V U +  b u  = 0, \\
(u,U)\equiv {\bf U},\, \bf{U}(0) = \bf{U}_0,\, \frac{d{\bf U}}{dt} (0) = {\bf U}'_0.
\end{array}
\end{equation}
can be constructed, and subsequently the  solutions of this  Cauchy  problem (\ref{Cauchy}) are obtained as  linear combinations
\[ {\bf U} (t) =  \sum_s U^s {\bf \Psi}_s \cos(\omega^b_s t + \varphi_s) = \Re  \sum_s U^s {\bf \Psi}_s  e^{i (\omega^b_s t + \varphi_s)} \]
of the  eigenmodes.
 But  we also can  reduce the above  homogeneous  equation (\ref{Cauchy}) to a pair  of  unperturbed  formally  inhomogeneous  equations,
 \begin{equation}
\label{Cauchy_unperturbed_cos_1}
m {u}_{tt} + v {u} + \sum_s U^s\,\,   b^+ \Psi_s\,\, \cos(\omega^b_s t + \varphi_s) \equiv  m {u}_{tt} + v {u} + {f},
\end{equation}
\begin{equation}
\label{Cauchy_unperturbed_cos_2}
M\, U_{tt} + V  U +  \sum_s  U^s\,\,  b\, \,\,\psi^b_s \,\, \cos(\omega^b_s t + \varphi_s) \equiv M\, U_{tt} + V  U  + {F} = 0,
\end{equation}
or a similar complex equation
\begin{equation}
\label{Cauchy_unperturbed_inhom_1}
m \vec{u}_{tt} + v \vec{u} + \sum_s U^s\,\,   b^+ \Psi_s\,\, e^{i(\omega^b_s t + \varphi_s)} \equiv  m \vec{u}_{tt} + v \vec{u} + \vec{f},
\end{equation}
\begin{equation}
\label{Cauchy_unperturbed_inhom_2}
M\, \vec{U}_{tt} + V  \vec{U} +  \sum_s  U^s\,\,  b\, \,\,\psi^b_s \,\, e^{i(\omega^b_s t + \varphi_s)} \equiv M\, \vec{U}_{tt} + V  \vec{U}  + \vec{F} = 0,
\end{equation}
with  inhomogeneities   $\vec{f} = b^+ U, \vec{F} =  b  \vec{u}$ obtained  via substitution,  for  $u, U$, the corresponding
 components of the  solution  ${\bf U}(t) =  \left(u, U \right)$ of the  Cauchy problem  (\ref{Cauchy}) for the  original oscillators system:
  \begin{equation}
\label{inhom_u}
{u}_{\vec{f}} (t) =  \sum_s   b^+\, \Psi_s \, {U}^s\,\,
\left[ \frac{ e^{i (\omega^b_s t + \varphi_s)}}{ m \lambda^b_s - v}\right],\,\,
{u}_f (t) =  \sum_s   b^+\, \Psi_s \, {U}^s\,\,
\left[ \frac{ \cos(\omega^b_s t + \varphi_s)}{ m \lambda^b_s - v}\right],\,\,
\end{equation}
  \begin{equation}
\label{inhom_U}
{U}_{\vec{F}} (t) =   \sum_s  b\, \psi_s \, U^s\,\,
\left[ \frac{ e^{i (\omega^b_s t + \varphi_s)}}{ M \lambda^b_s - V}\right],\,\,
{U}_{F} (t) =   \sum_s  b\, \psi_s \, U^s\,\,
\left[ \frac{ \cos(\omega^b_s t + \varphi_s)}{ M \lambda^b_s - V}\right]
\end{equation}

General solutions  of  the equations
 (\ref{Cauchy_unperturbed_cos_1}) -(
 \ref{Cauchy_unperturbed_inhom_2}) are  obtained  by adding  general solutions  of the corresponding  homogeneous equations
   $  m u_{tt} + v u = 0,\,\, M\, U_{tt} + V  U  = 0,\, $  to these  solutions  of the inhomogeneous equations.

The  last  formulas  allow  us to  calculate  partial values  of  energy of the ``small'' and  ``large'' oscillators  depending  on  time,  for instance.
\begin{equation}
\label{Energy}
E_f (U) = \frac{m |u_f'|^2  +  v | u_f|^2}{2},\,\, E_F (U) = \frac{M |U_F'|^2  +  V | U_F|^2}{2}
\end{equation}
under formally   ``exterior'' forces  $f, F $, defined by a linear combination  of the perturbed   normal modes  ${\bf U}$  of  the
perturbed  oscillator's system. 

The partial values  of   energy  of each oscillator  $u_f, U_F$ do not remain constant  in the course  of
 evolution, but depend on time, exposing beats while
the oscillators exchange energy  due to the  bond   $B =  {\rm antidiag} (b^+, b)$.  Beats   are  calculated using the  1D theory developed for  a periodic    exterior force as in  \cite[Section 22]{Landau_Mechanics}. 
To study  the  above (formally) complex version  of the  inhomogeneous  equation  (\ref{Cauchy_unperturbed_inhom_2}) we follow Landua \cite{Landau_Mechanics},
and  introduce the  data $\vec{\xi}_f = \sqrt{m} \,\,{u}'_{\vec{f}} + i \sqrt{v}\,\, u_{\vec{f}},\,\, {\xi}_f = \sqrt{m} \,\,{u}'_{{f}} + i \sqrt{v}\,\, u_{{f}} $
   and  re-write the  last equation
 (\ref{Cauchy_unperturbed_inhom_1}) as
\begin{equation}
\label{1D_inhom}
\sqrt{m}\,\, \vec{\xi}' - i \sqrt{v} \,\, \vec{\xi} =  \vec{f}(t),\,\, \sqrt{m}\,\, {\xi}' - i \sqrt{v} \,\, {\xi} =  {f}(t),\,\, \sqrt{v} =  \omega^0 \sqrt{m}.
\end{equation}
If  the  Cauchy problem   is solved, with an initial condition $\xi(0) = \xi_0$, then the energy (\ref{Energy}) of the small
oscillator for the real solution  ${\bf U}$ of the total problem,  is  calculated  as
\begin{equation}
\label{1D_inhom_xi}
 \bar{\xi}_f (t)\,\, \xi_f (t) = \overline{[\sqrt{m} {u}' + i \sqrt{v} {u}]}\,\,[ \sqrt{m} {u}' + i \sqrt{v} {u}] = 2E_f({u}).
\end{equation}

The last formula allows us to  estimate the  energy transfer  ``between the   modes  of unperturbed oscillators''. It  is sufficient
 to be able to  monitor the  time-depndent energy  of the small oscillator.
 
 \medskip

We now suggest  a  1D  version of the corresponding analysis, assuming that  $\mu = 1,\,\,K = C_1 \oplus C_1$. Back to
 discussion of  the  1D  oscillators, with  $\lambda^b_{1,2} =  \bar{\lambda} \pm h \equiv \lambda^b_{\pm}$, assume, that an approximate
 resonance  condition  is  satisfied,  
 \[ \omega^b_{1,2} =\omega^0_{1,2} + \frac{b^2}{4 m_1\,m_2\, \delta \lambda \omega^0_{1,2}}, \]
 and assume that
  the  solution of the  original  Cauchy  problem (\ref{Cauchy})
 is given as a linear combination of  the perturbed modes
 \[
 {\bf U} =
 \left(
 \begin{array}{c}
 u\\
 U
 \end{array}
 \right)=
 \sum_{s=1}^2 {\bf \Psi}^b_s U^s \cos (\omega^b_s t + \varphi_s).
 \]
Then the above inhomogeneous  evolution  equation for the small oscillator  is either
\[
m_1 u_1^{''} +  m (\omega^0)^2 u_1  + \sum_{s=1}^2 b^+ \Psi^b_s U^s \cos (\omega^b_s t + \varphi_s) = 0,
\]
or
\[
\sqrt{m} \left[\xi' - i \omega^0_1 \xi  \right] + \sum_{s=1}^2 b^+ \Psi^b_s U^s \cos (\omega^b_s t + \varphi_s) = 0,
\]
 and has  a partial solution
 \[
 u_1 = \sum_{s=1}^2 b^+ \Psi^b_s U^s \frac{\cos (\omega^b_s t + \varphi_s)}
 {m_1 \left[ (\omega^b_s)^2 - (\omega^0_1)^2 \right]}.
 \]
The corresponding  $\xi$ - function  is calculated  from the  equation
\[ \sqrt{m} \left[ \xi_1'  - i \xi_1 \omega_1^0\right] + \sum_{s=1}^2 b^+ \Psi^b_s U^s \cos (\omega^b_s t + \varphi_s) = 0,\]
 as
 \begin{eqnarray} \nonumber
 \xi_1 (t) &=& e^{i\omega^0_1\, t}
 \left[
 \xi_1 (0) - \frac{1}{\sqrt{m_1}}\int_0^t e^{-i\omega^0_1\,\,\tau }
 \sum_{s=1}^2 b^+ \Psi^b_s U^s \cos (\omega^b_s \tau + \varphi_s) d\tau  \right] \\ &\equiv & e^{i\omega^0_1\, t} \left[  \xi_1 (0) +  \hat{\xi}_1 (t) \right]
,  \label{xi-function}
 \end{eqnarray}
 with
\begin{eqnarray*}
 \xi_1 (0) & = &  \sqrt{m_1}\left[ {u'}_1 (0)  +  i \omega^0_1  u^0_{1} (0) \right] \\
& = & \frac{1}{\sqrt{m_1}}
 \sum_{s=1}^2 b^+ \Psi^b_s U^s \frac{ \left[- \omega^b_s \sin \varphi_s   + i
 \omega^0_1 \, \cos \varphi_s  \right]}
 { \left[ (\omega^b_s)^2 - (\omega^0_1)^2 \right]} \\
& = &  \frac{1}{\sqrt{m_1}}
 \sum_{s=1}^2 b^+ \Psi^b_s U^s \left [  \frac{\sin \varphi_s}{\omega_1^b + \omega_1^0}   + \frac{i \omega^0_1 e^{i\varphi_s}}{(\omega_s^b + \omega^0_1)(\omega_s^b - \omega^0_1)}                         \right]
\end{eqnarray*}
and  the  integral  $\hat{\xi}_1 (t) $ in  (\ref{xi-function})  is  calculated  as
\begin{eqnarray*}
\hat{\xi}_1 (t)  & = &  - \frac{1}{\sqrt{m_1}}\int_0^t e^{-i\omega^0_1\,\,\tau }
 \sum_{s=1}^2 b^+ \Psi^b_s U^s \cos (\omega^b_s \tau + \varphi_s) d\tau \nonumber \\
& = &  - \frac{1}{\sqrt{m_1}} \sum_{s = 1}^2\, b^+ \Psi^b_s U^s \; \cdot \\ &&
\left[ \frac{e^{i\varphi_s}}{2 i (\omega^b_s - \omega^0_1) } \left(e^{i (\omega^b_s - \omega^0_1)t} -
1 \right)  - \frac{e^{-i\varphi_s}}{2 i (\omega^b_s + \omega^0_1) } \left(e^{i (-\omega^b_s - \omega^0_1)t} - 1 \right)
\right] \\
& = &  - \frac{1}{\sqrt{m_1}} \sum_{s = 1}^2\, b^+ \Psi^b_s U^s\; \cdot \\ &&
\left[
e^{i\varphi_s} \, e^{i (\omega^b_s - \omega^0_1)\frac{t}{2}}\, \frac{\sin (\omega^b_s- \omega^0_1)\frac{t}{2}}{\omega^b_s- \omega^0_1} +
e^{-i\varphi_s} \, e^{-i (\omega^b_s + \omega^0_1)\frac{t}{2} }\, \frac{\sin (\omega^b_s + \omega^0_1)\frac{t}{2}}{\omega^b_s + \omega^0_1}
\right].
\nonumber \end{eqnarray*}

To compare the above  theoretical estimation  of time dependence  of the energy  of the ``small'' oscillator of  time,
 we  should  consider the  averaged energy  over a stepwise  system of  windows, such as  used  for  manufacturing the time-spectral cards discussed in Section 1.  Set
\begin{equation}
\label{average}
\overline{|\xi|^2} (T) = \frac{1}{\Delta}\int_{T - \Delta/2}^{T + \Delta/2} |\xi (t)|^2 dt,
\end{equation}
with  an  appropriate  choice  of  parameters for the averaging depending on  the basic  characteristics of $\lambda_{1,2}$  or  $\bar{\lambda}, \delta \lambda$.

\medskip

Let us first examine   the  time-dependence  of energy  of the ``small''oscillator  in  the simplest case  of two  1D  oscillators,  we  assume, that there exists  only  one
resonance eigenvalue of the perturbed  system, closest to the  unperturbed  eigenvalue of the  ``small'' oscillator. Further assuming that  $\frac{b^2}{m_1\,m_2} << (\delta \lambda)^2$,
we  may estimate the perturbed eigenvalues  and  eigenfrequencies   of the system as
\begin{eqnarray}
(\omega^b_1)^2 - (\omega^0_1)^2 &\approx& \bar{\lambda} +  \sqrt{(\delta \lambda)^2 + \frac{b^2}{m_1\,m_{2}}} - \lambda^0_1 =  \frac{b^2}{2 m_1\, m_2\, \delta \lambda\,\lambda^0_1 },
\nonumber \\
\label{eigen_values_perturbed}
(\omega^b_2)^2 - (\omega^0_1)^2 & \approx& \bar{\lambda} -  \sqrt{(\delta \lambda)^2 + \frac{b^2}{m_1\,m_{2}}} - \lambda^0_1 = - \delta \lambda - \frac{b^2}{2 m_1\, m_2 \,\delta \lambda\,\lambda^0_2 } \\
\nonumber \omega^b_1 - \omega^0_1 &\approx&  \frac{b^2}{4 m_1\, m_2\,\delta \lambda \,\omega^0_1 },
\,\,\,
 \omega^b_1 + \omega^0_1 \approx 2 \omega^0_1 + \frac{b^2}{4 m_1\, m_2\,\delta \lambda \,\omega^0_1 }
\\
\label{eigen_frequencies_perturbed}
\omega^b_2 - \omega^0_1 & \approx & \omega^0_2 - \omega^0_1 - \frac{b^2}{4 m_1\, m_2\,\delta \lambda \,\omega^0_2}\equiv \delta \omega  -  - \frac{b^2}{4 m_1\, m_2\,\delta \lambda \,\omega^0_2}.
\end{eqnarray}

Notice  the resonance terms with  small denominators  or/and   with slowly oscillating exponents  on the  window  give  the  crucial  contribution  to the average  (\ref{average}) while  the rapidly  oscillating  and  smooth  terms  may be neglected while  integrating over the window. Using the above representation for $\xi_1 (t) =  \xi_1 (0) + \hat{\xi} (t)$ we estimate the averaged  $\xi$-function ${\overline{|\xi|^2} (T) }$
\begin{eqnarray}\nonumber
{\overline{|\xi|^2} (T) } & = &  \frac{1}{\Delta}\int_{T - \Delta/2}^{T + \Delta/2} |\xi (t)|^2 dt \\ & = &  \frac{1}{m_1} \sum_{r,s = 1}^2 b^+ \Psi^b_r U^r \,\, b^+ \Psi^b_s U^s
\frac{1}{\Delta}   \int_{T - \Delta/2}^{T + \Delta/2} \overline{A} B  \end{eqnarray}
where 
\begin{eqnarray*}
A& = &    \frac{\sin \varphi_s}{\omega^b_s + \omega^0_1} + \frac{i \omega^0_1 e^{i\varphi_s}}{\lambda^b_s - \lambda^0_1}   + \left[  e^{i \varphi_s} e^{i (\omega^b_s -\omega^0_1)\frac{t}{2}}
  +  e^{-i \varphi_s} e^{- i (\omega^b_s +\omega^0_1)\frac{t}{2}} \right]
 \frac{\sin (\omega^b_s + \omega^0_1)\frac{t}{2}}{\omega^b_s + \omega^0_1}    \\
B & = & 
  \frac{\sin \varphi_r}{\omega^b_r + \omega^0_1} + \frac{i \omega^0_1 e^{i\varphi_r}}{\lambda^b_r - \lambda^0_1}   + \left[ e^{i \varphi_r} e^{i (\omega^b_r -\omega^0_1)\frac{t}{2}}
  +  e^{-i \varphi_r} e^{- i (\omega^b_r +\omega^0_1)\frac{t}{2}} \right]
 \frac{\sin (\omega^b_r + \omega^0_1)\frac{t}{2}}{\omega^b_r + \omega^0_1}   
\end{eqnarray*}
choosing the  window such that  the leading  resonance term  $\frac{\sin (\omega^b_r -\omega^0_1)\frac{t}{2}}{\omega^b_r -\omega^0_1}$ only  slightly deviates from  a constant on the window $|t-\tau|< \Delta$:

\[ \bigg|\sin (\omega^b_r -\omega^0_1)\frac{t}{2} -  \sin (\omega^b_r -\omega^0_1)\tau/2\bigg| \leq \Delta (\omega^b_r -\omega^0_1) \approx \frac{b^2\,\,\Delta}{4 m_1 m-2 \delta \lambda \omega^0_1}.\] On another hand, the  slowest  oscillation  of non-resonance  exponentials in the integrand is defined by the exponent 
\[  \frac{1}{\Delta}\int_{T- \Delta/2}^{T- \Delta/2} cos (\omega^b_s - \omega^0_1) t dt \approx
 \frac{1}{\Delta\,\, \delta \omega } \]
 and this should be small too.
\[
\frac{b^2\,\,\Delta}{4 m_1 m_2 \delta \lambda \omega^0_1} \approx  \frac{1}{\Delta\,\, \delta \omega }
\]
This now defines  the  width of the  optimal window giving the dependence  on  $\lambda^0_1, \lambda^0_2$ and the other parameters. 

Then, selecting this optimal window  and taking into account  only leading terms
of  the integrand,  we obtain an  approximate  estimation for the averaged  energy of the ``small'' oscillator  in dependence on time:
\begin{eqnarray}\nonumber
 \bar{|\xi|}(T) & = & \frac{1}{\Delta}\int_{T- \Delta/2}^{T+\Delta/2} |\xi|^2 (t) dt \approx |b^+ \Psi^s U^s|^2 dt \left[ \frac{|\omega^0_1|^2 }{\lambda^b_1 - \lambda^0_1} +
  \frac{\sin^2(\omega^b_1 -\omega^0_1)\frac{t}{2}}{(\omega^b_1 -\omega^0_1)^2} \right] \\
  & \approx&
 \label{averaged_energy}
 |b^+ \Psi^s U^s|^2 \left[  \frac{ 2 m_2 \delta \lambda \lambda^0_1}{b^2} +  \frac{\sin^2 (\frac{b^2}{4 m_1\,m_2 \delta\lambda \omega^0_1})}{ (\frac{b^2}{4 m_1\,m_2 \delta\lambda \omega^0_1})^2}  \right].
  \end{eqnarray}
\begin{figure}
[ht]
\begin{center}
 \includegraphics [width=4.5in] {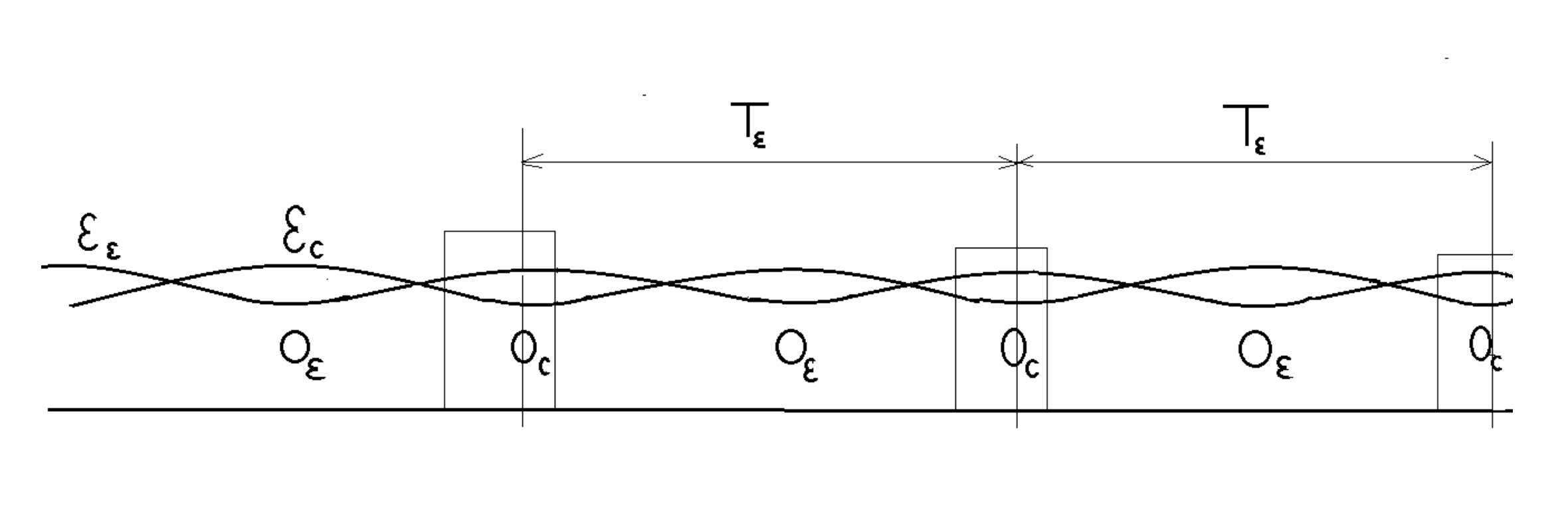}
\end{center}
\label{F:beating_symb}
\caption
{Symbolic diagram of the  energy content  ${\mathcal{E}}_{\varepsilon} (T),\,\,{\mathcal{E}}_{c}(T) $
of  the components  $u_{_{\varepsilon,c}} $  of the  perturbed dynamics $u$ with respect  to the
 de-localized  mode  $\Psi_c^0 $  on the complement $\Omega_c$ and  one of the localized  mode  on the
 active zone $\Omega_{\varepsilon}$.  Positions $O_{\varepsilon}, O_c$ of the  minima (and
  maxima) of the  energy content of the  localized and the delocalized  modes  alternate with
  opposite phases. The dangerous intervals of time  at the minima  $O_c$,  when  the destruction of $\Omega_{\varepsilon}$
   is  expected, are marked with  thin rectangles. }
\end{figure}
\vskip0.3cm

\section{ Appendix  1: Natural boundary conditions and the perturbed biharmonic wave equation. }
 In  our second section  we considered  a perturbed  biharmonic equation (\ref{perturbed_biharmonic}), borrowed  from \cite{Heisin_1967, Mikhlin_variational},
 and modelled the  normal stress  by  the boundary condition. Here we give more details describing  the  ``bridge'' connecting the model dynamics presented in
by us with the  corresponding chapter of classical mechanics of
  the plates and shells, see \cite{Timoshenko_circular,Timoshenko_plates_shells,Mikhlin_variational} and also \cite{Greenwood_1989,KP_1986}.  We keep in mind, that the linear  theory of small
  oscillations  of tectonic plates will only shed light on an initial  phase  of the  process  which may  lead to the   catastrophic  results. But we hope that
    the  the mechanical  realization of the preliminary  small oscillation   model  of the  resonance  process may be able to  preview some initial  features  of  the  catastrophic
  phase  of the  process. In reality,  constructing     the mechanical  model of the resonance interaction of  the  SGO  modes of tectonic plates  is a problem to resolve  based on experiment,
     with  use  of far more detailed
  mechanical details, see  for instance \cite{Timoshenko_plates_shells, Reddy_2007, Maurer}.  We  only attempt  here   a  first step  in  this direction,  by considering
    two   thin  tectonic plates  developing  typical  kinds  of stresses when colliding  under    fluctuation  of the
 rotation speed  of Earth  an/or  the  convection flow  in  the liquid underlay (astenosphere).  To further simplify  our analysis, we   assume  that  both  $\Omega_{\varepsilon},\Omega_c$  are  Kirchhoff plates, see \cite{Timoshenko_plates_shells}, and endure
  different kinds of stresses.
 We  assume that  the  normal pressure  and  the corresponding bending  dominate  the  potential energy  of the large plate $\Omega_{\varepsilon}$ and  the
 tangential (shearing)  component
  of the  stress dominates the  potential energy  of the small plate. The  shearing  part   of the stress  is defined  by   the  tension $T $ in  the middle
   plane  of
   with the components $T_x, T_y, T_{xy}$, satisfying the equilibrium conditions, see Mikhlin's book, \cite[Chapter 4, Section 28]{Mikhlin_variational} which we use as a basic reference in what comes.  The equilibrium conditions are 
   \begin{equation}
  \label{equilibrium}
  \frac{\partial T_x}{\partial x } +  \frac{\partial T_{xy}}{\partial y } = 0,\,\,
  \frac{\partial T_y}{\partial y } +  \frac{\partial T_{xy}}{\partial x } = 0,
  \end{equation}
and the corresponding   quadratic form
\[  \langle w, T w \rangle =  \int_{\Omega_{\varepsilon}}
  \left[ w_x T_x w_x  + 2 w_x  T_{xy} w_y  +  w_y T_y w_y  \right] d\Omega \]
  defines, with  the kinematic boundary condition  $w\bigg|_{\Gamma} = 0$,  a symmetric operator $T$ of second order on smooth  functions  $w$.
  The operator $T$ is  positive  for  stretching tension and  negative  for  compressing tension, see
  \cite{Heisin_1967}.

  This may  cause an instability  of the tectonic plate, see for instance the  analysis  of a numerical example for an elliptic plate under a
   compressive tension in \cite[Chapter 8, Section 72]{Mikhlin_variational}.  The  materials composing
  tectonic plates  can't resist the stretching tension, so  hereafter we assume that  the tension $T$ is  compressing, and the corresponding  second
   order operator $T$ is negative.

 We  begin with considering of  the non-interacing plates  $\Omega_{\varepsilon}, \Omega_c$ under normal  bending stress $F$, and  the compressing  stress $T$,
 assuming that each of them is
  elastically fixed    at the  common boundary  $\Gamma_{\varepsilon}$, see below.  The  zero boundary condition is applied on all boundaries, free reclining
  conditions are applied  on  the remote part
  of boundary  of  $\Omega_c$, and  the  elastic fixture on both sides  $\Gamma_{\varepsilon,c}$  of the common boundary  $\Gamma$ is modeled by  the corresponding
   functionals
   $\int_{\Gamma} \beta_{\varepsilon, c} \left(\frac{\partial w_{\varepsilon, c}}{\partial n}\right)^2  d\Gamma$.

      The  equilibrium
   deformations $w_{\varepsilon}, w_c$ are  found  by  minimizing  of
  the  corresponding quadratic functionals on  the subspace of the  virtual strains,
  subject to the kinematic boundary condition $ w\bigg|_{\Gamma_{\varepsilon}} = 0 $
\begin{eqnarray}
 W (u)  \nonumber  &=&  D \int_{\Omega} |\Delta w|^2 d\Omega -  2(1-\sigma)  D \int_{\Omega}\left({w}_{xx} w_{yy}  - |u_{xy}|^2\right) d\Omega \\ && +
\int_{\Gamma} \beta \left(\frac{\partial w}{\partial n}\right)^2  d\Gamma + H \langle w, T w \rangle - 2 \int_{\Omega} {w} F  d\Omega. \label{potential energy}
\end{eqnarray}

Here $\sigma$ is  the Poisson coefficient  $0<\sigma<1$, $D$ is the flexural ridigity,  $D = \frac{H^2 E}{12 (1-\sigma^2)}$,$E$ is the Young modulus, $H $
is the thickness
of the  plate and  $\beta> 0 $ is the parameter defining the elastic  contact  of the plate with environment. The  second integral  of the  Monge-Ampere  form
 (the curvature  of  the  surface  $z = w(x,y)$) can be  presented  as an integral on the boundary, see \cite{Birman_1956}, with regard of the above kinematic
  boundary condition $w\bigg|_{\Gamma} = 0$:
\[
2 \int_{\Omega}\left({w}_{xx} w_{yy}   - |w_{xy}|^2\right) d\Omega =  \int_{\Gamma}\frac{\left(\frac{\partial w}{\partial n}\right)^2}{r (\gamma)} d\Gamma
\]

For  the circular $\Gamma_{\varepsilon}$ we  assume $r(\gamma) = \varepsilon$ and  $H \langle w, T w \rangle =  H^{-2} D  \Delta w$. Calculation of the first
variation of the energy functional yields
 the Euler equation
\begin{equation}
\label{Euler}
D \Delta^2 w  + T w   =  F
\end{equation}
with the kinematic and  the natural (under the above elastic $\beta$- bound)   condition on the each side of the common boundary, e.g.
 \begin{equation}
D \Delta w  +  \beta \frac{\partial w}{\partial n} - D\frac{1-\sigma}{r}  \frac{\partial w}{\partial n} \bigg|_{\Gamma_{\varepsilon}} = 0.
\label{elastic_bound}
\end{equation}

To derive an equation  for  the small ocillation  we  consider, for each plate, the   Lagrangian, associated with  the thin plate $\Omega$,
subject to  the  above elastic bound and  the kinematic boundary conditions $u\bigg|_{\Gamma} = 0$
\begin{eqnarray}
\mathcal{L} & = &   \int_0^t\int_{\Omega}  D|\Delta w|^2 d\Omega\, dt  -
   \int_0^t\int_{\Omega}\left[
  2D (1- \sigma)\left(w_{xx}w_{yy} - w_{xy}^2\right) -
  \rho H w^2_t\right] d\Omega\, dt \nonumber \\  
  && +
 2 \int_0^t\int_{\Gamma} \beta \left(\frac{\partial w}{\partial n}\right)^2 d\Gamma dt +  \int_0^T  \langle w, T w \rangle dt,
 \end{eqnarray}
we obtain  the perturbed  biharmonic  wave equation  as   the Euler equation for the  critical  points of the Lagrangian:

\begin{equation} 
\rho H w_{tt} +  D  \Delta^2 w  +  T w  = 0,
\end{equation}
with the   ``$\beta $- natural''  free reclining boundary conditions on both sidesof the common boundary
\begin{equation}
w\bigg|_{\Gamma} = 0,\,\,
 D \Delta w + \beta \frac{\partial w}{\partial n} - \frac{D(1-\sigma)}{r} \frac{\partial w}{\partial n} \bigg|_{\Gamma} = 0,
\label{bc_free}
\end{equation}
and free reclining boundary condition on the complementary part of the boundary  of $\Omega_{c}$.

While  considering  a circular plate $\Omega:  0< r < a$, assume that the active zone $\Omega_{\varepsilon}$ is centered in  $\Omega$, so that  the complement
$\Omega_{c}$  is a ring $0< r < a$, as described earlier.  Generally we have the  Lagrangian
\begin{eqnarray}
\mathcal{L} & = &  \int_0^t\int_{\Omega}  D|\Delta w|^2 d\Omega\, dt  +  \int_0^t \int_{\Gamma_{\varepsilon}}
 \beta \left(\frac{\partial w}{\partial n}\right)^2 d\Gamma_{\varepsilon} dt
  - \int_0^t\int_{\Omega}
  \rho H  \left(\frac{\partial w}{\partial t}\right)^2 d\Omega\, dt \nonumber \\
&&   
- 2D (1- \sigma) \int_0^t \left(w_{xx}w_{yy} - w_{xy}^2\right) d\Omega dt  + \int_{0}^t \langle w, T w \rangle   dt,
 \end{eqnarray}
This gives  the perturbed  biharmonic  wave equation  as the  Euler equation for the  critical  points  of the Lagrangian  on $\Omega_c$.

\begin{equation}
\label{Euler_biharmonic}
\rho H w_{tt} +  D  \Delta^2 w  +  T w  = 0,
\end{equation}
and
 \begin{equation}
\label{Euler_biharmonic_varepsilon}
\rho H w_{tt} +  D  \Delta^2 w  +  T w  = 0,
\end{equation}
on $\Omega_{\varepsilon}$, with $T <0$
and the boundary condition  with regard of the  elastic bond:

\[
w\bigg|_{\Gamma_{\varepsilon}} = 0,\,\,  \beta \frac{\partial w}{\partial n} D\Delta w  \bigg|_{\Gamma_{\varepsilon}}  -
\frac{D (1-\sigma)}{\varepsilon} \frac{\partial w}{\partial n} \bigg|_{\Gamma_{\varepsilon}} = 0.
\]
 For the  complement we  neglect the  tangential compression, but  keep  the bending and the  elastic bond, so that the  Lagrangian  is reduced to
\begin{equation}
 \label{Lagrangian_biharmonic}
\mathcal{L}_{c} =  \int_0^T\int_{\Omega}  D|\Delta w|^2 d\Omega\, dt  +  \int_0^T \int_{\Gamma_{\varepsilon}} \beta
 \left(\frac{\partial w}{\partial n}\right)^2 d\Gamma_{\varepsilon} dt
  - \int_0^T\int_{\Omega}
  \rho H  \left(\frac{\partial w}{\partial t}\right)^2 d\Omega\, dt,
 \end{equation}
we obtain  the perturbed  biharmonic  wave equation  as  Euler equation for the  critical  points  of the Lagrangian   $L_c$:

\begin{equation} 
\rho H w_{tt} +  D  \Delta^2 w   = 0,
\end{equation}
with the  boundary condition
\[
w\bigg|_{\Gamma_c} = 0,\,\,  \beta \frac{\partial w}{\partial n} - D \frac{1-\sigma}{r} \frac{\partial w}{\partial n} + D\Delta w  \bigg|_{\Gamma_{\varepsilon}} = 0.
\]

We  make  now one more  simplifying assumption,
assuming that the  operator $L_{\varepsilon}$  on   the  active zone is defined on  circular harmonics of zero order  $ n=0$   (independent of the angular variable), and
 the  operator $L_c$ on the complement is  defined on the linear span  of  circular harmonics of the first order,  $e^c_1 =  \pi^{-1}\cos \varphi,\,\,\,e^s_1 =\pi^{-1}\sin \varphi$,
 so that the corresponding   eigenfunctions  are  spanned  by  the circular harmonics.

  The  unperturbed  spectral problems, associated with the  1D differential equations
 \[
 L_c^{c,s} w =  D_{c} \Delta^2 w  + Tw  =  \omega^2 H_c \rho w
  \]
 have 4  linearly  independent solutions   $J_1, H^1_1\equiv H_1, I_1,K_1$  with factors  $\sin \varphi, \cos \varphi$.
 
 We denote them by   $J^{c,s}_1,  H^{c,s}_1, I^{c,s}_1,K^{c,s}_1$.

Selecting  appropriate  elastic bonds, we  may set the parameters $\beta_0, \beta_1^c,\beta_1^s$ so that  the  unperturbed spectral problems   for $L_{\varepsilon}
 \oplus L^c_c \oplus L_c^s$
has a multiple  eigenfrequency  $\nu = (2\pi)^{-1} \omega$
for $\delta = 0$.   Then the  perturbed  spectral problems,  associated with the unperturbed  operator
\[
\left[
\begin{array}{ccc}
D \Delta_0^2  +  Q \Delta &  0 & 0 \\
0 &  D \Delta_1^2  &  0 \\
0 &  0 &  D \Delta_1^2
\end{array}
\right]
\equiv L^0_{\varepsilon} \oplus  L^c_{1}\oplus L^s_{1}
\]
and separate   $\beta$-natural  boundary conditions  on the common part of the boundary.   For the perturbed  problem, defined  by the same differential expression
\begin{equation}
\label{spectral_perturbed_CS}
\left[
\begin{array}{ccc}
\beta_0 - D_{\varepsilon}\frac{1-\sigma}{r_{\varepsilon}} &  \kappa^c e_0 \rangle \langle e^c_1 &  \kappa^s e_0 \rangle \langle e^s_1 \\
 \kappa^c e^c_1 \rangle \langle e_0  & \beta_1^c - D_c\frac{1-\sigma}{r_{\varepsilon}} &  0 \\
  \kappa^s e^s_1 \rangle \langle e_0  & 0 & \beta_1^s - D_c \frac{1-\sigma}{r_{\varepsilon}}
\end{array}
\right]
\left[
\begin{array}{c}
\frac{\partial \Psi_{0}}{\partial n}\\
\frac{\partial \Psi
^c_{1}}{\partial n}\\
\frac{\partial \Psi^s_{1}}{\partial n}
\end{array}
\right]
=
\Delta
\left[
\begin{array}{c}
 \Psi_{0}\\
 \Psi^c_{1}\\
 \Psi^s_{1}
\end{array}
\right]
\end{equation}

Here  $\kappa^{c,s}$ are   small real  parameters. The multiple eigenfrequency is split  into  the  starlet of simple perturbed eigenfrequencies   with
eigenfunctions  constructed  as  linear combinations
 $\Psi_0$ of  $J_0, I_0$  on $\Omega_{\varepsilon}$, see Section 4, and a linear combination $\Psi^{c,s}_c$ of  the
 $J^{c,s}_1,  H^{c,s}_1, I^{c,s}_1,K^{c,s}_1$ with coefficients
 found from the above analog (\ref{spectral_perturbed_CS}) of the  free reclining boundary conditions and  the unperturbed
 boundary condition  on the  complementary part   of  the  boundary.

Contact:  g.j.martin@massey.ac.nz

\end{document}